\documentclass[preprint, 5p, times]{elsarticle}
\biboptions{sort&compress}
\usepackage{nomencl}

\usepackage{microtype} 
\usepackage{xurl} 
\usepackage{etoolbox}
\renewcommand\nomgroup[1]{%
  \item[\bfseries
    \ifstrequal{#1}{A}{Sets}{%
    \ifstrequal{#1}{B}{Indices}{%
    \ifstrequal{#1}{C}{Parameters}{%
    \ifstrequal{#1}{D}{Decision Variables}{%
    \ifstrequal{#1}{E}{Auxiliary Variables}{Other}}}}}%
  ]\vspace{10 pt}
}

\newcommand{\nomunit}[1]{%
\renewcommand{\nomentryend}{\hspace*{\fill}#1}}
\makenomenclature

\usepackage{hyperref}
\hypersetup{
    colorlinks=true,
    linkcolor=blue,
    filecolor=magenta,      
    urlcolor=cyan,
    pdfpagemode=FullScreen,
    }

\usepackage{amsmath}
\allowdisplaybreaks 
\usepackage{afterpage}
\usepackage{amssymb}
\usepackage{booktabs}
\usepackage{tabularx}
\usepackage{adjustbox} 
\usepackage{threeparttable}
\usepackage{siunitx}
\usepackage{chemformula}
\usepackage{caption}
\usepackage{subcaption}
\usepackage{comment}

\usepackage{perpage}
\MakePerPage{footnote}
\usepackage[official]{eurosym}

\usepackage{algorithm2e}

\SetKwFor{ParFor}{parfor}{do}{end}%

\usepackage{algpseudocode}

\RestyleAlgo{ruled}
\SetKw{GoTo}{go to}

\usepackage{xspace}
\newcommand{\Gas}{\textbf{G}}
\newcommand{\PV}{\textbf{P}}
\newcommand{\Wind}{\textbf{W}}
\newcommand{\Hydro}{\textbf{H}}
\newcommand{\Nuclear}{\textbf{N}}

\newcommand{\Inverter}{\textbf{I}}
\newcommand{\Battery}{\textbf{B}}

\newcommand{\Scenarios}{\mathcal{S}}
\newcommand{\ScenariosSimulation}{\mathcal{S}'}
\newcommand{\Regions}{\mathcal{R}}
\newcommand{\Plants}{\mathcal{P'}}
\newcommand{\Technologies}{\mathcal{P}}
\newcommand{\TimeSteps}{\mathcal{T}}
\newcommand{\AllHours}{\mathcal{H}}
\newcommand{\ValidLinks}{\mathcal{L}}

\newcommand{\MaxCapacity}[2]{\text{cap}_{#1,#2}^{\text{max}}}
\newcommand{\cf}[4]{\text{cf}_{#1,#2,#3,#4}}
\newcommand{\timeStep}{\text{ts}}
\newcommand{\plannedOutageRate}{\beta}
\newcommand{\hydroReservoir}[1]{\text{rl}_{#1}^{\text{max}}}
\newcommand{\resInflow}[2]{\text{ri}_{#1,#2}}
\newcommand{\distance}[2]{\text{dist}_{#1,#2}}
\newcommand{\load}[2]{\text{load}_{#1,#2}}
\newcommand{\transmissionEfficiency}[2]{\eta_{#1, #2}^{\text{T}}}
\newcommand{\efficiency}[1]{\eta_{#1}}
\newcommand{\emissionFactor}[1]{\text{ef}_{#1}}
\newcommand{\dischargeTime}{\text{dt}}
\newcommand{\fuelCost}[1]{\text{c}_{#1}^{\text{FL}}}
\newcommand{\varCost}[1]{\text{c}_{#1}^{\text{V}}}
\newcommand{\fixedCost}[1]{\text{c}_{#1}^{\text{F}}}
\newcommand{\invCost}[1]{\text{c}_{#1}^{\text{I}}}
\newcommand{\crf}[1]{\text{crf}_{#1}}
\newcommand{\crfTransmission}{\text{crf}^{\text{T}}}
\newcommand{\transmissionCost}{\text{c}^{\text{T}}}
\newcommand{\carbonTax}{\text{c}^\text{tax}}
\newcommand{\sheddingCost}{\text{c}^{\text{S}}}
\newcommand{\loadLossCost}{\text{BM}}
\newcommand{\sheddingCapRatio}{\text{sr}}
\newcommand{\prob}[1]{\pi_{#1}}
\newcommand{\SC}{\text{SC}}
\newcommand{\LoL}{\text{LoL}}
\newcommand{\lossPerKm}{\alpha^{\text{Loss}}}

\newcommand{\discountRate}{\text{i}}
\newcommand{\lifeTime}{\text{T}}

\newcommand{\OperationalCost}[1]{\text{OC}_{#1}}
\newcommand{\FixedCost}{\text{FC}}
\newcommand{\TotalLoadShedding}[1]{\text{TS}_{#1}}

\newcommand{\TotalEmission}[1]{\text{TE}_{#1}}
\newcommand{\InvestmentCost}{\text{IC}}
\newcommand{\Electricity}[4]{\text{E}_{#1,#2,#3,#4}}
\newcommand{\ElectricityShedding}[3]{\text{ES}_{#1,#2,#3}}
\newcommand{\ElectricityLoadLoss}[3]{\text{EL}_{#1,#2,#3}}
\newcommand{\Capacity}[2]{\text{Cap}_{#1,#2}}

\newcommand{\TransmissionCapacity}[2]{\text{Cap}_{#1,#2}^{\text{T}}}
\newcommand{\BatteryLevel}[3]{\text{BL}_{#1,#2,#3}}
\newcommand{\ReservoirLevel}[3]{\text{RL}_{#1,#2,#3}}
\newcommand{\Transmission}[4]{\text{T}_{#1,#2,#3,#4}}
\newcommand{\PlannedOutage}[3]{\text{PO}_{#1,#2,#3}}
\newcommand{\UnPlannedOutage}[3]{\text{UO}_{#1,#2,#3}}

\newcommand{\oneNuclearPlantSize}{\text{m}}
\newcommand{\Month}{\mathcal{M}}
\newcommand{\Sizes}{\mathcal{N}}
\newcommand{\noMonthlySamples}{\text{NM}}
\newcommand{\noNuclearPlants}{\text{n}}
\newcommand{\noSamplesforPercentiles}{\text{N}^{\text{P}}}
\newcommand{\noSamplesinRobustModel}{\text{N}^{\text{R}}}
\newcommand{\noSimulations}{\text{N}}

\newcommand{\confidenceLevel}{\alpha}
\newcommand{\SelectedSample}{\text{z}}
\newcommand{\SelectedSampleBold}{\textbf{z}}
\newcommand{\SelectedSampleBar}{\bar{\text{z}}}
\newcommand{\SelectedSampleBarBold}{\bar{\textbf{z}}}
\newcommand{\outageSample}[2]{\text{os}_{#1, #2}}

\newcommand{\annualOutagePercentile}{\text{AOP}}
\newcommand{\maxNoOutagePercentile}{\text{MOP}}

\newcommand{\startHour}{\text{h}^{s}}
\newcommand{\sampleIndex}{\text{id}}

\newcommand{\outageTime}{\text{ot}}

\newcommand{\currentOutageStartHour}{\bar{\text{h}}^{s}}
\newcommand{\hoursPerMonth}[1]{\text{h}^{\text{M}}_{#1}}
\newcommand{\montlyOutageSamples}[2]{\text{os}^{\text{M}}_{#1, #2}}
\newcommand{\outageSamples}[2]{\text{os}_{#1, #2}}
\newcommand{\outageSamplesBold}{\textbf{os}}
\newcommand{\simulatedOutageSamples}[3]{\text{os}_{#1, #2, #3}^{'}}
\newcommand{\simulatedOutageSamplesBold}{\outageSamplesBold^{'}}

\newcommand{\CapacityBar}[2]{\overline{\text{Cap}}_{#1,#2}}
\newcommand{\CapacityBarBold}{\overline{\textbf{Cap}}}

\newcommand{\Card}[1]{\lvert #1 \rvert}

\DeclareSIUnit{\TonCOtwo}{\text{t}_{\ch{CO2}}}
\DeclareSIUnit{\MTonCOtwo}{\text{Mt}_{\ch{CO2}}}
\DeclareSIUnit{\KgCOtwo}{\text{kg}_{\ch{CO2}}}
\DeclareSIUnit{\Currency}{\text{\euro}\xspace}
\DeclareSIUnit{\year}{\text{y}}
\DeclareSIUnit{\MW}{\text{MW}}
\DeclareSIUnit{\GW}{\text{GW}}
\DeclareSIUnit{\KW}{\text{kW}}
\DeclareSIUnit{\KWh}{\text{kW}}
\DeclareSIUnit{\MWh}{\text{MWh}}
\DeclareSIUnit{\GWh}{\text{GWh}}
\DeclareSIUnit{\TWh}{\text{TWh}}
\DeclareSIUnit{\MCurrency}{\text{M\euro}\xspace}

\newcommand{\Favorable}{\text{F}\xspace}
\newcommand{\Average}{\text{A}\xspace}
\newcommand{\Unfavorable}{\text{U}\xspace}

\newcommand{\Deterministic}{D\xspace}
\newcommand{\Stochastic}{WU\xspace}
\newcommand{\Robust}{WNU\xspace}

\newcommand{\D}{\textit{Deterministic}\xspace}
\newcommand{\WU}{\textit{Weather Uncertainty}\xspace}
\newcommand{\WNU}{\textit{Weather and Nuclear Uncertainty}\xspace}

\newcommand{\NoOutageSimulation}{\textit{No-outage}\xspace}
\newcommand{\NormalSimulation}{\textit{Normal}\xspace}
\newcommand{\UnfavorableSimulation}{\textit{Unfavorable Weather}\xspace}
\newcommand{\DunkelflauteSimulation}{\textit{Dunkelflaute}\xspace}

\newcommand{\SEN}{$\text{SE}^\text{N}$\xspace}
\newcommand{\SES}{$\text{SE}^\text{S}$\xspace}
\newcommand{\DEN}{$\text{DE}^\text{N}$\xspace}
\newcommand{\DES}{$\text{DE}^\text{S}$\xspace}
\newcommand{\DK}{\text{DK}\xspace}
\newcommand{\BNL}{\text{BNL}\xspace}
\newcommand{\PL}{\text{PL}\xspace}

\begin{document}

\begin{frontmatter}



\title{Modeling Robust Energy Systems Considering Weather Uncertainty and Nuclear Power Failures: A Case Study in Northern Europe}


\affiliation[chalmers]{
            organization={Department of Space, Earth and Environment},
            addressline={Chalmers University of Technology}, 
            postcode={412 96}, 
            state={Göteborg},
            country={Sweden}}
            
\author[chalmers]{Kamran Forghani\corref{cor1}}

\author[chalmers]{Xiaoming Kan}
\author[chalmers]{Lina Reichenberg}
\author[chalmers]{Fredrik Hedenus}

\cortext[cor1]{Corresponding author}
\ead{kamranf@chalmers.se}


\begin{abstract}
Capacity expansion models used for policy support have increasingly represented both the variability and uncertainty of weather-dependent generation (wind and solar). However, although also uncertain, as demonstrated by the performance of the French nuclear power fleet in 2022, uncertainty arising from nuclear power outages has been largely neglected in the literature. This paper presents the first capacity expansion model that considers uncertainty in nuclear power availability caused by unplanned outages. We propose a mathematical model that combines a scenario-based stochastic optimization approach (to deal with weather-related uncertainties) with a data-driven adjustable robust optimization approach (to deal with nuclear failure-related uncertainties). The robust model represents the bulky behavior of nuclear power plants, with large (1~\unit{\GW}) units that are either on or off, while at the same time letting the model decide on the optimal amount of nuclear capacity. We tested the model in a case for Northern Europe (seven nodes) with a time resolution of 1250 time steps. Our findings show that nuclear power outages do, in fact, impose a vulnerability on the energy system if not considered in the planning phase. Our proposed model performs well and finds solutions that prevent Loss-of-Load (at a price of robustness of 0.6\%), even in more extreme weather conditions. Robust solutions are characterized by a higher capacity of gas plants, but, perhaps surprisingly, nuclear power capacity is barely affected.

\end{abstract}


\begin{keyword}
    weather uncertainty \sep nuclear outage \sep stochastic optimization \sep robust optimization \sep electricity systems\sep  capacity expansion
\end{keyword}

\fntext[fn1]{This work has been submitted to \emph{Applied Energy} for possible publication.}
\end{frontmatter}
\section{Introduction} \label{sec:Introduction}

Renewable energy technologies, such as wind and solar power, along with nuclear power, are potentially key energy sources in low-carbon energy systems. However, debates persist within the scientific community, as well as in public and political spheres, regarding how their variability limits the extent to which energy systems may rely on them, see \citep{kan2020cost, sepulveda2018role}. While the uncertainty in wind and solar power outputs is often highlighted, the risk of unplanned outages of nuclear power plants receives comparatively little attention. Even less explored is the potential interaction between these two uncertainties. Despite their significance, incorporating uncertainty into large-scale energy system planning is not a common practice, even though its importance has been highlighted by previous studies \citep{larsson1993developing, decarolis2017formalizing}. In this study, we develop a novel method that combines stochastic and robust optimization to address uncertainties in wind, solar, and nuclear power outputs for robust energy system design.

The need to account for these uncertainties has become increasingly evident. The recent energy crisis in Europe vividly highlights the urgent need for robust energy systems that can reliably meet demand even under unusual circumstances \citep{mivsik2022eu, brodny2023assessing}. However, the energy system analyses that inform policymakers' decisions rely heavily on deterministic approaches, which often result in solutions that are not resilient to even regular operational uncertainties (e.g., inter-annual weather variability) and may underestimate the true cost of a robust energy system. 

To address uncertainty in energy system planning, several studies have explored the importance of weather year selection \citep{gotske2024designing, collins2018impacts} and the impact of prolonged periods of low solar and wind output \citep{ruhnau2022storage, holtinger2019impact}. These studies typically adopt a deterministic approach and rely on scenario or sensitivity analysis \citep{hsia1994formal}, which helps assess the effects of different assumptions about uncertain inputs. Similar methods include global sensitivity analysis \citep{iooss2015review} and Monte Carlo simulations \citep{raychaudhuri2008introduction}. However, deterministic approaches do not explicitly incorporate uncertainty into the decision-making process \citep{wallace2000decision, yue2018review}.

To address decision-making under uncertainty, stochastic and robust optimization methods are increasingly being used to investigate various aspects of the energy system, such as investment costs \citep{fu2021effects, dominguez2020planning, moret2020decision, bistline2013electric}, fuel costs \citep{moret2020decision}, discount rates \citep{moret2020decision}, and energy demand \citep{dominguez2020planning, dimanchev2024consequences}. More relevant to our focus, some studies have incorporated wind and solar output uncertainty into the design of future energy systems. \citet{seljom2015short} made early efforts to incorporate wind power output uncertainty in an energy system model for Denmark. \citet{perera2020quantifying} developed a stochastic-robust optimization method to account for wind and demand uncertainties in Sweden. \citet{seljom2021stochastic}  incorporated stochastic wind and solar outputs into Norway’s long-term energy model, finding that stochastic optimization leads to 13\% less solar and 2\% more wind capacity, with a system cost difference of less than 1\%. \citet{verastegui2019adaptive} applied adaptive robust optimization to address operational uncertainties in renewable power outputs and electricity demand in Chile, showing that a robust design increases system costs by 8\%–20\% compared to a deterministic approach. Based on the literature, weather uncertainty's impact on system costs seems to vary significantly across different models and contexts.

Regarding nuclear power outages, \citet{murphy2020resource} examined the impact of temperature-dependent forced outages of thermal power generators on the operation of a regional electricity market in the US. They found that accounting for temperature-dependent failures increases the required reserve capacity from 16\% to 23\%. In contrast to temperature-dependent outages in warm climate zones, nuclear power outages caused by technical failures are more prevalent across all countries and far less predictable. These failures can result in full outages lasting several months, significantly affecting the reliability of the electricity system \citep{iaea2023operating}. For instance, in France, a record 26 of its 56 reactors were offline due to pipe corrosion during the energy crisis in Europe in 2022 \citep{nyt2022france}. However, to the best of our knowledge, no studies have yet incorporated these uncertainties into decision-making for energy system planning.

In this paper, we present an energy system model for capacity expansion that incorporates two types of uncertainties: weather-related uncertainty arising from the variability in wind and solar power outputs due to changing weather conditions and failure-related uncertainty arising from nuclear power plant outages. We use a scenario-based approach to address weather-related uncertainties, while a robust optimization method is applied to deal with failure-related uncertainty. To assess the impact of these uncertainties on energy system planning, we conduct experiments using the developed model for the northern European electricity system.

The aims of the paper are to:
\begin{enumerate} 
    \item Develop a methodology to account for both weather uncertainties and unplanned nuclear power outages in energy system planning.   
    \item Quantify the direction and approximate magnitude of impacts on cost and energy system design from incorporating weather uncertainties and unplanned nuclear power outages into long-term energy system planning. 
    \item Assess the \emph{cost} of resilience, as well as the \emph{value} of designing for resilience if something does not go as planned. 
\end{enumerate}    

The remainder of this paper is organized as follows: In Section~\ref{sec:LiteratureReview}, we review typical methods used to address uncertainties related to component failure in power systems. Section~\ref{sec:ProblemDescription} explains the base energy system model and key assumptions and formulates the problem as a scenario-based two-stage program. We then derive the robust counterpart by incorporating unplanned nuclear failure uncertainty and present the methodology for solving the robust model. Section~\ref{sec:Experimentaldesign} introduces the case study, weather scenarios, and simulation cases, followed by presenting and discussing the results and comparing system performance across different cases in terms of system cost, capacity mix, and loss of load. Finally, Section~\ref{sec:Conclusion} provides conclusions and directions for future research.
\section{Review of studies on component failure in energy systems} \label{sec:LiteratureReview}

In this paper, we study a capacity expansion problem for electricity generation facing two types of uncertainties: uncertainty in the wind and solar power outputs (weather-related uncertainty) and uncertainty in the availability of nuclear power plants (uncertainty related to component failure in energy systems). Scenario-based stochastic optimization approaches have widely been used in the literature to account for weather-related uncertainties, see, e.g., \citep{backe2021stable, bagheri2022stochastic, ma2024optimal}. Since a relatively large and accurate dataset of wind and solar capacity factors is accessible, this approach can effectively represent different weather-year conditions using scenarios. However, limited attention has been paid to component failure, more specifically nuclear power plant failure, in the literature. The main focus of this section is to review the existing studies that utilize stochastic and robust optimization approaches to deal with component failures. For more comprehensive reviews on stochastic and robust optimization applications in electricity systems, readers are referred to \citep{zakaria2020uncertainty, qiu2022application, rahim2022overview}.

In real-world electricity systems, due to unpredictable events, components such as transmission lines, power generation units, and equipment face the risk of failure. Researchers have utilized different approaches to represent the uncertain operating states of these components when failures happen; examples of such studies can be found in \citep{khosrojerdi2016method, alabdulwahab2015stochastic}.

Although scenario-based optimization approaches are straightforward to implement, they have limitations in effectively accounting for component failure in power systems. Section~\ref{sec:RobustModel} discusses these limitations in more detail. Aside from using scenarios, another approach to representing component failures is to model them using box uncertainty. For example, \citet{ratanakuakangwan2021hybrid} assumed that available dispatchable capacity is influenced by the risk of power plant outages and applied \citet{soyster1973convex}'s robust optimization framework (based on box uncertainty) to ensure the total dispatchable capacity meets the projected peak demand. Despite its simplicity, the conventional box uncertainty approach is highly conservative, resulting in the maximum price of robustness; see \citep{bertsimas2004price, asgari2024data} for more details. 

In many cases, optimization problems are single-period, where decisions are made all at once, and there is no opportunity to adapt. Static Robust Optimization (RO) models are well-suited in those cases. These models focus on minimizing risk based on the worst-case realization of uncertain parameters. In contrast, for multi-period optimization problems, which happens to be more common in electricity systems, \citet{ben2004adjustable} proposed the idea of Adjustable Robust Optimization (ARO), also known as two-stage or adaptive robust optimization. In ARO, decision-making occurs in two stages. The first stage involves \textit{Here-and-Now} variables (e.g., investment decisions on power plants, transmission lines, and storage facilities), which are determined before uncertainty is revealed. The second stage involves \textit{Wait-and-See} variables (e.g., power generation schedule and load shedding), which are decided after the uncertainty is resolved.

\citet{wang2013two} proposed an ARO approach to address the contingency-constrained unit commitment problem, considering contingencies in both generators and transmission lines. To obtain a robust unit commitment schedule, the objective was to minimize the total generation cost in a multi-bus power grid under an $N\text{--}k$ security criterion. The $N\text{--}k$ security criterion specifies that for a system with $N$ elements, there are $\binom{N}{k}$ possible contingencies, and no more than $k$ failures can occur simultaneously. In a similar study conducted by \citet{zhang2019two}, researchers dealt with this issue by extending the $N\text{--}k$ security criterion to also include the probabilities of component failures in the problem. These studies mainly concentrated on the modeling of uncertainty sets considering the simultaneous failure of multiple components in the system. More extensive analysis was made by, for instance, \citet{guo2018islanding}, where not only singular events were considered but also uncertainties related to the duration of component failures. They investigated energy management in microgrids under uncertainties in both the power outputs of renewable energy sources and the duration of microgrid islanding event periods \footnote{events where the microgrid is disconnected from the network due to interruption in transmission lines}. They assumed that the islanding time is an uncertain parameter and the microgrid undergoes the islanding mode for $k$ out of $N$ periods. Similarly, \citet{gholami2019proactive} studied microgrid islanding events while considering uncertainties in both the islanding time and duration.

ARO approaches, such as the ones reviewed earlier, offer a more accurate representation of system functionality under uncertain conditions. By employing ARO, it is possible to develop more resilient energy systems, often leading to less conservative solutions compared to static RO. This approach is therefore used in this study to address nuclear power plant outages.

\subsection{Research gaps in nuclear outage modeling} \label{sec:research_gap}

In capacity expansion problems for electricity systems, the power outputs of the generation units (such as wind, solar, and nuclear) could vary due to weather conditions and technical failures. Researchers have already used stochastic and robust optimization approaches to deal with uncertainty in wind and solar power outputs. However, to the best of our knowledge, the uncertainty in nuclear power output due to unplanned outages has not been investigated. A summary of the reviewed papers addressing component failures in electricity systems is provided in Table~\ref{tab:literature_summery}. As this table shows, there is only one study, \citep{ratanakuakangwan2021hybrid}, dedicated to component failures in capacity expansion problems where the authors used a conventional robust optimization approach (box uncertainty) to deal with uncertainty in the capacity of dispatchable generating units. Modeling uncertainty in the availability of nuclear power plants requires considering both the failure start time and duration (referred to as temporal uncertainty in the table). Additionally, multiple plants may fail simultaneously (denoted as spatial uncertainty in the table), and historical data show that failure rates vary across different months. Moreover, the number of nuclear plants in the system is a decision variable that must be determined by solving the capacity expansion model. However, existing models developed to deal with component failure models do not account for these factors simultaneously. To address this research gap, this paper introduces an ARO approach to represent uncertainty in nuclear power outages. 

\begin{table*}[htbp]
    \caption{Overview of studies addressing uncertainties related to component failures in power systems.} \label{tab:literature_summery}
    \resizebox{\textwidth}{!}{
        \begin{threeparttable}
            \begin{tabular}{llcccccccc} \toprule
                 &  &  & \multicolumn{4}{c}{Uncertainty related to component failures} & \multicolumn{2}{c}{Other uncertainties} &  \\ \cmidrule(lr){4-7} \cmidrule(lr){8-9}
                Paper & \begin{tabular}[c]{@{}c@{}}Problem \\ scope\end{tabular} & \begin{tabular}[c]{@{}c@{}}Investment \\ decisions\end{tabular} & \begin{tabular}[c]{@{}c@{}}Uncertain \\ parameters\end{tabular} & \begin{tabular}[c]{@{}c@{}}Spatial \\ uncertainty\end{tabular} & \begin{tabular}[c]{@{}c@{}}Temporal \\ uncertainty\end{tabular} & \begin{tabular}[c]{@{}c@{}}Uncertainty \\ modeling method\end{tabular} & \begin{tabular}[c]{@{}c@{}}Uncertain \\ Parameters\end{tabular} & \begin{tabular}[c]{@{}c@{}}Uncertainty \\ modeling method\end{tabular} & \begin{tabular}[c]{@{}c@{}}Solution \\ approach\end{tabular} \\ \midrule
                \cite{khosrojerdi2016method} & Power supply chain & \checkmark & \begin{tabular}[c]{@{}c@{}}Generator and \\ transmission line\\ outages\end{tabular} & \checkmark &  & Scenario-based RO\tnote{b} & -- & -- & MIP solver \\ \midrule
                \cite{alabdulwahab2015stochastic} & \begin{tabular}[c]{@{}c@{}}Unit commitment with\\ natural gas network\end{tabular} &  & \begin{tabular}[c]{@{}c@{}}Generator and \\ transmission line\\ outages\end{tabular} & \checkmark & \checkmark & \begin{tabular}[c]{@{}c@{}}Scenario-based\\ SO\end{tabular} & Load & Scenario-based SO & MIP solver \\ \midrule
                \cite{ratanakuakangwan2021hybrid} & Capacity expansion & \checkmark & Reliable capacity &  &  & RO\tnote{c} & \begin{tabular}[c]{@{}c@{}}Capacity factors\\ Load\\ Capital and fixed \\ O\&M cost\end{tabular} & Scenario-based RO\tnote{b} & MIP solver \\ \midrule
                \cite{wang2013two} & Unit commitment &  & \begin{tabular}[c]{@{}c@{}}Generator and \\ transmission line\\ outages\end{tabular} & \checkmark &  & \begin{tabular}[c]{@{}c@{}}ARO with $N\text{--}k$\\  security criterion\end{tabular} & -- & -- & \begin{tabular}[c]{@{}c@{}}Bender’s\\ decompositions\end{tabular} \\ \midrule
                \cite{zhang2019two} & Unit commitment &  & Generator outages & \checkmark & \checkmark & \begin{tabular}[c]{@{}c@{}}ARO with $N\text{--}k$ \\ and $\alpha$-cut\\ security criterion\end{tabular} & \begin{tabular}[c]{@{}c@{}}Capacity factors\\ Load\end{tabular} & \begin{tabular}[c]{@{}c@{}}ARO with box-type \\ uncertainty set\end{tabular} & \begin{tabular}[c]{@{}c@{}}Bender’s\\ decomposition\end{tabular} \\ \midrule
                \cite{guo2018islanding} & Microgrid scheduling &  & \begin{tabular}[c]{@{}c@{}}Micro grid \\ islanding time\end{tabular} &  & \checkmark & ARO\tnote{a} & Capacity factors & ARO\tnote{a}  & \begin{tabular}[c]{@{}c@{}}Column and constraint\\ constraint generation\end{tabular} \\ \midrule
                \cite{gholami2019proactive} & Microgrid scheduling &  & \begin{tabular}[c]{@{}c@{}}Micro grid \\ islanding time\\  and duration\end{tabular} &  & \checkmark & ARO\tnote{a} & \begin{tabular}[c]{@{}c@{}}Capacity factors\\ Load\\ Market prices\end{tabular} & ARO\tnote{a}  & \begin{tabular}[c]{@{}c@{}}Column and constraint\\ constraint generation\end{tabular} \\ \midrule
                This paper & Capacity expansion & \checkmark & \begin{tabular}[c]{@{}c@{}}Nuclear power\\ outages\end{tabular} & \checkmark & \checkmark & \begin{tabular}[c]{@{}c@{}}Data-driven ARO\\ with $N\text{--}k$ security\\ criterion\end{tabular} & Capacity factors & Scenario-based SO & \begin{tabular}[c]{@{}c@{}}Heuristic \\ decomposition\end{tabular} \\ \bottomrule
            \end{tabular}
        \begin{tablenotes}
            \item[a,b,c] \citet{bertsimas2004price} approach [a]; \citet{mulvey1995robust} approach [b]; \citet{soyster1973convex} approach [c]
        \end{tablenotes}
    \end{threeparttable}
    }
\end{table*}

\section{Problem description and mathematical model} \label{sec:ProblemDescription}

The capacity expansion model used in this study is a modified version of the model proposed in \cite{kan2021impacts}. The model represents investment in power generation capacities, batteries, and transmission lines. Power generation sources include photovoltaic solar (\PV), wind (\Wind), hydropower (\Hydro), and nuclear (\Nuclear), along with a natural gas turbine (\Gas). The power produced using these technologies can be stored in batteries for future use. Each technology is characterized by parameters such as investment cost, operational cost, \ch{CO2} emissions factor, and capacity factor. Transmission lines can also be established to transfer the electricity among a set of regions $\Regions$.

\begin{figure}[htbp]
    \centering
    \includegraphics[width=\columnwidth]{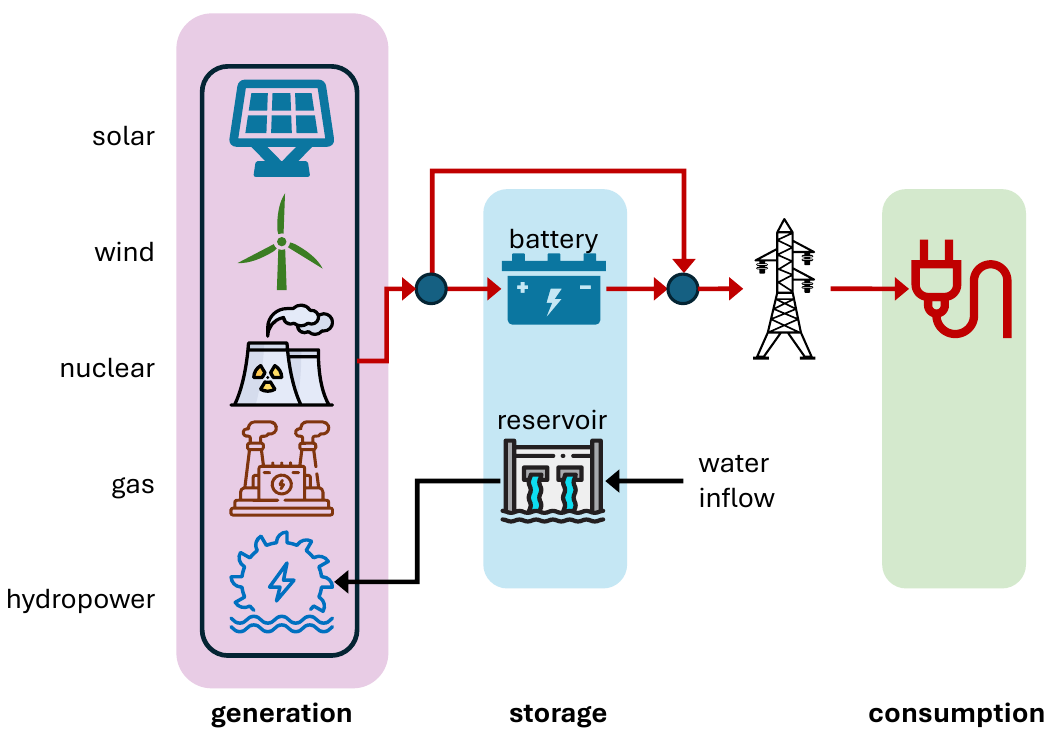}
    \caption{Schematic of interactions in the proposed energy systems.}
    \label{fig:Electricity_System_Schematic}
\end{figure}

Since the electricity load and the inputs for wind, solar, and hydropower technologies are variable in time, in our model, we use annual time series data to capture the variations in those inputs. Having an hourly temporal resolution makes the model computationally heavy to solve. Therefore, we use a parameter called the time step, $\timeStep$, to adjust the model's temporal resolution. The objective of the models is to find a capacity mix that minimizes the levelized total cost of the system. Load shedding is allowed in our model, and the shed load is penalized in the objective function by a parameter called the load shedding cost, $\sheddingCost$. The amount of shed load at each time step, however, must not exceed a certain fraction, $\sheddingCapRatio$, of the electricity load. There is no constraint on \ch{CO2} emission levels, and a carbon tax approach is applied in the objective function to penalize emissions. The electricity system and the interactions among its components are depicted in Figure~\ref{fig:Electricity_System_Schematic}. 

In real-world power systems, nuclear power plants may become unavailable due to technical failures and maintenance requirements. To model nuclear plant availability, we categorize nuclear power outages into two types: \emph{planned outages}, resulting from scheduled maintenance, and \emph{unplanned outages}, mainly caused by unexpected technical failures. Since planned outages can be scheduled in advance, we introduce the parameter $\plannedOutageRate$, which represents the fraction of time that needs to be allocated for scheduled nuclear plant maintenance. This parameter is directly incorporated into the stochastic model discussed in Section~\ref{sec:StochasticModel}. In contrast, unplanned outages are unpredictable and beyond operational control. The consequences of such outages may become significantly problematic, especially if a substantial portion of nuclear capacity becomes unavailable during periods of low wind and solar generation. To address this, Section~\ref{sec:RobustModel} extends the stochastic model by integrating unplanned nuclear outages through a data-driven ARO approach.

It should be mentioned that the notation list provided in \ref{sec:appendix_notation} is used throughout the paper in the development of the mathematical models.

\subsection{\WU model}\label{sec:StochasticModel}

We use a scenario-based stochastic optimization approach to capture the uncertainty in wind and solar power outputs (i.e., weather-related uncertainties), a method widely used in the literature (see, for example, \citep{backe2021stable, bagheri2022stochastic, ma2024optimal}). Since a relatively large and accurate dataset of wind and solar capacity factors is available, we can select several weather years as scenarios representing different weather-year conditions. Some of these scenarios can be chosen in a way that realistically captures extreme weather events based on historical data. The mathematical model for the proposed capacity expansion problem, which incorporates uncertainty in wind and solar capacity factors using scenarios, is given in \ref{sec:stochastic_model}. This model is referred to as the \WU (\Stochastic) model, and later, in Section~\ref{sec:RobustModel}, it is extended to incorporate the uncertainty in nuclear power availability. It should be noted that, in Section~\ref{sec:Results}, the term \D (\Deterministic) model is used to refer to a specific case of the \Stochastic model. In this case, a single weather year close to the average weather conditions is used as a scenario in the \Stochastic model for making comparisons between stochastic and deterministic approaches.

\subsection{\WNU model}\label{sec:RobustModel}

In this section, the \Stochastic model described in Section~\ref{sec:StochasticModel} is extended to address unplanned nuclear power outages. Scenario-based approaches, such as those presented in \citep{alabdulwahab2015stochastic, khosrojerdi2016method}, have limitations in representing the uncertainty associated with component failures due to the following reasons:

\begin{itemize} 
    \item \textbf{Complexity of the scenario tree}: Since the number and layout of components (e.g., nuclear power plants) are decision variables, an upper bound must be set, and distinct scenarios must be created to represent their state over time (since nuclear power plants are treated as discrete units\footnote{Unlike the simultaneous reduction in wind and solar power output, nuclear failures typically do not occur at the same time, which adds complexity in scenario generation for nuclear plants.}). This approach is impractical, as the number of scenarios grows exponentially, making stochastic optimization computationally intractable. 
    \item \textbf{Underestimation of consequences}: Stochastic optimization may underestimate the risks of rare but extreme events, such as simultaneous failures of multiple components during periods of low wind and solar output, which could lead to a loss of load. In contrast, robust optimization can provide preparedness against worst-case scenarios, thus mitigating the risks of loss of load. 
\end{itemize}

Additionally, approaches such as those proposed in \citep{wang2013two, zhang2019two} (which employ the $N\text{--}k$ security criterion) and \citep{guo2018islanding, gholami2019proactive} (focused on the time and duration of outages) are not directly applicable for modeling nuclear power outages. The reason is that, firstly, in a capacity expansion problem, the number and layout of nuclear plants are decision variables rather than predefined inputs. Secondly, when modeling nuclear power outage uncertainty, factors such as outage start time, duration, the number of simultaneous outages, and their time correlation must be considered\footnote{Historical data shows that outages are influenced by weather conditions, e.g., higher failure rates during hot seasons due to cooling requirements.}. These factors make existing uncertainty modeling methods (see Table~\ref{tab:literature_summery}) less effective for representing nuclear power outages.

Therefore, we employ an ARO approach to address nuclear power outage uncertainty. In this respect, a data-driven uncertainty set is combined with a discrete-state uncertainty set to define the uncertainty in the ARO model. A data-driven uncertainty set uses empirical data to construct a polyhedron that models uncertainty in parameters. These sets are used to obtain robust solutions to optimization problems by capturing possible variations in uncertain parameters while relying on observed historical data. This helps reduce the conservativeness of robust solutions; see \citep{zhang2021data, ordoudis2021energy, shen2022data, lu2024distributionally}.

We utilize historical nuclear power outage data \citep{iaea2023operating} to simulate unplanned nuclear power outages throughout the year. Let $\noMonthlySamples$ denote the number of outage samples available in the historical data, and let $\montlyOutageSamples{i}{m}$ denote the unplanned outage time (in hours) for the $i$-th sample in month $m$. We can use these samples to simulate unplanned nuclear power outages on an hourly basis, denoted by $\outageSample{i}{t}$. To do so, we apply Algorithm~\ref{func:GenerateOutageSamples}($\noSimulations$) to generate $\noSimulations$ samples, representing whether a nuclear plant is offline (represented by 1) or active (represented by 0) at each hour of the year. It should be noted that to convert monthly outages, $\montlyOutageSamples{i}{m}$, to hourly samples, Algorithm~\ref{func:GenerateOutageSamples}($\noSimulations$) assumes that each monthly unplanned outage involves one incident that starts randomly within $[\startHour, \startHour + \hoursPerMonth{j} - \montlyOutageSamples{i}{j}]$, where $\startHour$ and $\hoursPerMonth{j}$ represent the starting hour and the total hours in the corresponding month, respectively.

\begin{function}[htbp]
    \caption{GenerateOutageSamples(\noSimulations)} \label{func:GenerateOutageSamples}
    \DontPrintSemicolon
    \SetArgSty{textrm}
    $\outageSample{i}{t} \gets 0, \quad \forall i \in \left\{1, 2, \ldots, \noSimulations \right\}, \; t \in \TimeSteps$\;
    \For{$i \gets 1$ \KwTo $\noSimulations$}{
        $\startHour \gets 1$\;
        $\sampleIndex \gets \textbf{Rand}(1, \noMonthlySamples)$\;
        \For{$j \in \Month$}{
            $\outageTime \gets \montlyOutageSamples{\sampleIndex}{j}$\;
            \If{$\outageTime > 0$}{
                $\currentOutageStartHour \gets \startHour + \textbf{Rand}(0, \hoursPerMonth{j} - \outageTime)$\;
                \For{$h \gets \currentOutageStartHour$ \KwTo $\currentOutageStartHour + \outageTime - 1$}{
                    $\outageSample{i}{t} \gets 1$\;
                }
            }
            $\startHour \gets \hoursPerMonth{j} + \startHour$\;
        }
    }
    \Return $\outageSamplesBold$\;
\end{function}

The amount of unplanned nuclear power outages depends on the number of nuclear plants in the system. Let’s assume that the number of nuclear plants in the system, i.e., $\sum_{r \in \Regions} \Capacity{r}{\Nuclear}$, is equal to $\noNuclearPlants \in \Sizes$, where $\Sizes$ defines the set of possible numbers of nuclear plants in the system. To adjust the conservativeness of the model, we introduce a parameter called the confidence level, $\confidenceLevel$, and impose two sets of constraints as uncertainty budgets in the uncertainty set. One constraint ensures that the annual unplanned nuclear power outages remain below a threshold, denoted by $\annualOutagePercentile$. The other constraint ensures that no more than $\maxNoOutagePercentile$ outages occur simultaneously at any hour. For any given number of nuclear plants, $\noNuclearPlants$, and confidence level $\confidenceLevel \in (0, 1)$, we use Algorithm~\ref{func:Percentiles}($\noNuclearPlants$, $\confidenceLevel$, $\noSamplesforPercentiles$) to perform simulations and determine $\annualOutagePercentile$ and $\maxNoOutagePercentile$. In this algorithm, the parameter $\noSamplesforPercentiles$ represents the number of outage samples generated to estimate $\annualOutagePercentile$ and $\maxNoOutagePercentile$. With these explanations, we define the following data-driven discrete-state uncertainty set, denoted by $\mathcal{U}(\noNuclearPlants, \confidenceLevel, \noSamplesforPercentiles)$, for the ARO model.

\begin{align}
& \mathcal{U}(\noNuclearPlants, \confidenceLevel, \noSamplesforPercentiles) := \Bigg\{ 
    \SelectedSample_{i, s} \in \{0, 1\}: \nonumber \\
    & (\annualOutagePercentile, \maxNoOutagePercentile) = \protect\textbf{\ref{func:Percentiles}}(\noNuclearPlants, \confidenceLevel, \noSamplesforPercentiles), \label{eq:percentiles} \\
    & \sum_{i=1}^{\noSamplesinRobustModel} \SelectedSample_{i, s} = \noNuclearPlants, \label{eq:no_nuclear_plants} \\
    & \sum_{i=1}^{\noSamplesinRobustModel} \sum_{t \in \AllHours}\outageSample{i}{t} \cdot \SelectedSample_{i, s} \leq \annualOutagePercentile, \label{eq:annual_outages} \\
    & \sum_{i=1}^{\noSamplesinRobustModel} \outageSample{i}{t} \cdot \SelectedSample_{i, s} \leq \maxNoOutagePercentile, \quad \forall i \in \left\{1, 2, \ldots, \noSamplesinRobustModel \right\}, \; s \in \Scenarios, \; t \in \AllHours \Bigg\} \label{eq:maximum_outages}
\end{align}

The so-called uncertainty set comprises three sets of constraints. Constraint~\eqref{eq:percentiles} sets the annual and maximum number of unplanned nuclear power outages by performing $\noSamplesforPercentiles$ simulations via Algorithm~\ref{func:Percentiles}(.)\footnote{$\noSamplesforPercentiles$ needs to be sufficiently large to ensure the accurate calculation of $\annualOutagePercentile$ and $\maxNoOutagePercentile$. For our experiments in Section~\ref{sec:Experimentaldesign}, we set $\noSamplesforPercentiles = 100,000$.}. Constraint~\eqref{eq:no_nuclear_plants} ensures that, out of $\noSamplesinRobustModel$ outage samples, one is selected for each nuclear plant in every scenario. Essentially, this constraint is inspired by the data-driven uncertainty set and determines both the outage time and duration for each nuclear plant\footnote{$\noSamplesinRobustModel$ also needs to be sufficiently large to ensure that the uncertainty set adequately covers all potential outage scenarios. This parameter, however, affects the computational complexity of the model. For our experiments in Section~\ref{sec:Experimentaldesign}, we set $\noSamplesinRobustModel = 5,000$.}. Constraints~\eqref{eq:annual_outages} and~\eqref{eq:maximum_outages} serve as two uncertainty budgets that respectively limit the annual and hourly outages.

\begin{function}[htbp]
    \caption{Percentiles($\noNuclearPlants$, $\confidenceLevel$, $\noSimulations$)} 
    \label{func:Percentiles}
    \DontPrintSemicolon
    \SetArgSty{textrm}
    $\textbf{O}^{\max} \gets \emptyset$ and $\textbf{O}^{\text{T}} \gets \emptyset$ \;
    \For{$i \gets 1$ \KwTo $\noSimulations$}{
        $\outageSamplesBold \gets \protect\textbf{\ref{func:GenerateOutageSamples}}(\noNuclearPlants)$ \;
        $\textbf{O}^{\max} \gets \textbf{O}^{\max} \cup \left\{\max_{t \in \AllHours}{\sum_{j = 1}^{\noNuclearPlants}\outageSamples{j}{t}}\right\}$ \;
        $\textbf{O}^{\text{T}} \gets \textbf{O}^{\text{T}} \cup \left\{\sum_{t \in \AllHours}{\sum_{j = 1}^{\noNuclearPlants}\outageSamples{j}{t}}\right\}$ \;
    }
    $\maxNoOutagePercentile \gets \textbf{Quantile}(\textbf{O}^{\max}, \confidenceLevel)$ \;
    $\annualOutagePercentile \gets \textbf{Quantile}(\textbf{O}^{\text{T}}, \confidenceLevel)$ \;
    \Return $\maxNoOutagePercentile, \annualOutagePercentile$ \;
\end{function}

Now, assuming that the number of nuclear plants in the system, $\noNuclearPlants$, belongs to the set $\Sizes$, the Min-Max-Min robust model of the problem, under the uncertainty set $\mathcal{U}(\noNuclearPlants, \confidenceLevel, \noSamplesforPercentiles)$, can be expressed as follows. Henceforth, this model will be referred to as the \WNU (\Robust) model.

\begin{align}
    & \hypertarget{RobustProblem}{\text{\Robust~model:}}  \nonumber \\
    & \SC(\confidenceLevel, \noSamplesforPercentiles, \noSamplesinRobustModel)  = \min_{\noNuclearPlants \in \Sizes}  \max_{\SelectedSampleBold \in \mathcal{U}(\noNuclearPlants, \confidenceLevel, \noSamplesforPercentiles)} \min \Bigg\{ \InvestmentCost + \FixedCost  \nonumber \\
    & \qquad + \sum_{s \in \Scenarios} \prob{s} \left(\OperationalCost{s} + \carbonTax \cdot \TotalEmission{s} + \sheddingCost \cdot \TotalLoadShedding{s}\right) \label{eq:objective_function_robust}\\
    & \text{Subject to: (\ref{eq:max_capacity})--(\ref{eq:electricity_storage_balance}), (\ref{eq:gen_capacity_wind_PV})--(\ref{eq:non_negativity_4})}  \nonumber \\
    & \Electricity{s}{r}{\Nuclear}{t} \leq \Capacity{r}{\Nuclear} \cdot \timeStep - \PlannedOutage{s}{r}{t} - \UnPlannedOutage{s}{r}{t}, \quad \forall s \in \Scenarios, \; r \in \Regions, \; t \in \TimeSteps \label{eq:gen_capacity_nuclear_robust} \\
    & \sum_{r \in \Regions}{\Capacity{r}{\Nuclear}} = \oneNuclearPlantSize \cdot \noNuclearPlants  \label{eq:total_nuclear_capacity} \\
    & \sum_{r \in \Regions}{\UnPlannedOutage{s}{r}{t}} = \oneNuclearPlantSize \cdot \timeStep \sum_{i=1}^{\noSamplesinRobustModel} \outageSample{i}{t} \cdot \SelectedSample_{i, s}, \quad \forall s \in \Scenarios, \; t \in \TimeSteps \label{eq:unplanned_load_loss_robust} \\
    & \UnPlannedOutage{s}{r}{t} \geq 0, \quad \forall s \in \Scenarios, \; r \in \Regions, \; t \in \TimeSteps \Bigg\} \label{eq:non_negativity_5}
\end{align}

In the \Robust model, the objective function in Equation~\eqref{eq:objective_function_robust} minimizes the worst-case cost of the system under the uncertainty set $\mathcal{U}(\noNuclearPlants, \confidenceLevel, \noSamplesforPercentiles)$ and confidence level $\confidenceLevel$. The innermost optimization problem in Equation~\eqref{eq:objective_function_robust} minimizes the $\SC$, given that the number of nuclear plants in the system is $\noNuclearPlants$, and $\SelectedSampleBold$ are the selected samples. All of the constraints in the innermost optimization problem are the same as those in the stochastic model presented in Section~\ref{sec:StochasticModel}, except for Constraint~\eqref{eq:gen_capacity_nuclear}, which needs to be replaced with Constraint~\eqref{eq:gen_capacity_nuclear_robust} to account for unplanned nuclear power outages as well. Three additional constraints must be considered. Constraint~\eqref{eq:total_nuclear_capacity} ensures that the total capacity of the nuclear plants in the system equals the number of nuclear plants multiplied by the unit capacity of a nuclear plant, $\oneNuclearPlantSize$, in \unit{\GW}. Constraint~\eqref{eq:unplanned_load_loss_robust} ensures that the total unplanned nuclear power loss (in \unit{\MWh}) in each time step is determined according to the selected outage samples. Finally, constraint~\eqref{eq:non_negativity_5} states that the decision variable corresponding to the unplanned nuclear power outages is non-negative.

\subsection{Simulation model} \label{sec:SimulationModel}

We perform simulations to evaluate the quality of solutions derived from solving the models. This allows us to assess the solutions more accurately using a relatively large historical dataset. For a given solution, $\CapacityBarBold$ (representing the capacity mix of various technologies) and $\noNuclearPlants$ (the number of nuclear plants in the corresponding capacity mix), the simulation involves independently solving the \hyperlink{SimulationProblem}{Simulation problem} (given in \ref{sec:simulation_model}) for each year $s \in \ScenariosSimulation$, where $\ScenariosSimulation$ is the set of simulation years. Since some solutions may fail to meet the demand constraint in certain simulation years, we introduce the decision variable $\ElectricityLoadLoss{s}{r}{t}$ to quantify the load loss in simulation year $s$, at region $r$, in time step $t$. This load loss is then penalized using a sufficiently large cost, $\loadLossCost$, in the objective function. After solving the optimization problem for each $s \in \ScenariosSimulation$, the average values of the system costs and loss-of-loads are calculated to represent the simulation outcome.

\subsection{Solution method} \label{sec:SolutionMethod}

Solving the \Robust model presented in Section~\ref{sec:RobustModel} to optimality is computationally challenging. This is because the proposed mathematical model is a nonlinear optimization (due to the tri-level Min-Max-Min structure of the objective function in \eqref{eq:objective_function_robust}), making it difficult for current solvers to solve the problem effectively. Although this model can be linearized using duality theory\footnote{The model can be converted into a mixed-integer program by writing the dual of the innermost problem as a maximization problem. Let $\lambda_{s, h}$ denote the dual variable associated with Constraint~\eqref{eq:unplanned_load_loss_robust}. From duality theory, we would have $\oneNuclearPlantSize \cdot \timeStep \sum_{s \in \Scenarios}\sum_{t \in \TimeSteps}\sum_{i=1}^{\noSamplesinRobustModel} \outageSample{i}{t} \cdot \SelectedSample_{i, s} \cdot \lambda_{s, h}$ in the objective function. The bi-linear term $\outageSample{i}{t} \cdot \SelectedSample_{i, s}$ can be linearized using an auxiliary variable and additional constraints.}, see also \citep{wang2014robust, qiu2022application}, solving the reformulated model still remains challenging, even when the capacity mix, $\Capacity{r}{p}$, is predetermined. One could argue that to reduce the number of binary variables, we can reduce the number of samples in the \Robust model, i.e., $\noSamplesinRobustModel$, or increase the time step parameter, i.e., $\timeStep$, as it significantly affects the number of variables and constraints in the problem. This, however, makes the model unrealistic to the extent that the resulting solutions become far from reality.

Thus, to solve the \Robust model more effectively, we propose a heuristic approach that mitigates these complications by solving it iteratively through a master problem and a sub-problem. When the number of nuclear plants, $\noNuclearPlants$, and the binary variable, $\SelectedSampleBold$, in the \Robust model are fixed and denoted as $\SelectedSampleBarBold$, the model simplifies to a linear program (LP), referred to as the \hyperlink{Masterproblem}{Master problem}($\noNuclearPlants$, $\SelectedSampleBarBold$).

\begin{align}
    & \hypertarget{MasterProblem}{\text{Master problem}}(\noNuclearPlants, \SelectedSampleBarBold)\text{:} \min \SC \\ 
    & \text{subject to: (\ref{eq:max_capacity})--(\ref{eq:electricity_storage_balance}), (\ref{eq:gen_capacity_wind_PV})--(\ref{eq:non_negativity_4}), (\ref{eq:gen_capacity_nuclear_robust}), (\ref{eq:total_nuclear_capacity}), (\ref{eq:non_negativity_5})}  \nonumber \\
    & \qquad \sum_{r \in \Regions}{\UnPlannedOutage{s}{r}{t}} = \oneNuclearPlantSize \cdot \timeStep \sum_{i=1}^{\noSamplesinRobustModel} \outageSample{i}{t} \cdot \SelectedSampleBar_{i, s}, \quad \forall s \in \Scenarios, \; t \in \TimeSteps \label{eq:unplanned_load_loss_master_problem} 
\end{align}

Assume a known capacity mix, denoted by $\CapacityBarBold$, consisting of $\noNuclearPlants$ nuclear plants. For this partial solution, to identify the extreme periods when nuclear power outages maximize the system cost, we can formulate a Max-Min optimization problem for each scenario $s$. However, this results in the same computational complexity issue discussed earlier due to the nonlinearity of the model. Our heuristic simplifies this process by avoiding the direct solution of the Max-Min problem. Instead, a simplified problem is solved to identify potential periods when nuclear power outages are more likely to result in loss of load. Given a known capacity mix $\CapacityBarBold$, we substitute the Max-Min objective function with another function that correlates with the highest system cost. A suitable alternative is assigning nuclear power outages to those time steps when the net electricity load, $\gamma_{t}$, as given by Equation~\eqref{eq:net_load}, is at its maximum. Our preliminary tests revealed that a normalized version of this parameter, shown in Equation~\eqref{eq:normalized_net_load}, yields better results (i.e., generally resulting in a higher system cost). 

\begin{align}
    &\gamma_{t} = \max \Bigg\{0, \sum_{r \in \Regions} \load{r}{t} - \sum_{p \in \left\{\Hydro, \Gas \right\}}\CapacityBar{r}{p} - \sum_{p \in \left\{\Wind, \PV \right\} }\ \crf{p} \cdot \CapacityBar{r}{p} \nonumber \\ 
    & \qquad - \oneNuclearPlantSize \cdot (\noNuclearPlants - \maxNoOutagePercentile) \Bigg\}, \quad \forall t \in \AllHours  \label{eq:net_load}\\
    & \bar{\gamma}_{t} = \frac{\gamma_{t}}{\max_{t \in \AllHours} \left\{1, \gamma_{t} \right\} }, \quad \forall t \in \AllHours \label{eq:normalized_net_load}   
\end{align}

After calculating the parameter $\bar{\gamma}_{t}$ using Equation~\eqref{eq:normalized_net_load}, for each scenario $s$ from the \Robust model, we derive a MIP model given by \hyperlink{SubProblem}{Sub-problem}(.) to identify the potential problematic periods. The objective function of this model maximizes a weighted outage time, where the weights correspond to the net load. In other words, this model prioritizes assigning nuclear power outages to the time steps with the highest net load.

\begin{align}
    & \hypertarget{SubProblem}{\text{Sub-problem}}(\noNuclearPlants, \annualOutagePercentile, \maxNoOutagePercentile, \noSamplesinRobustModel, \CapacityBarBold, s)\text{:}  \nonumber \\ 
    & \max \sum_{i=1}^{\noSamplesinRobustModel} \sum_{t \in \TimeSteps} \gamma_{t} \cdot \outageSample{i}{t} \cdot \SelectedSample_{i} \label{eq:objective_function_sub_problem}\\
    & \text{subject to:} \nonumber \\
    & \quad \sum_{i = 1}^{\noSamplesinRobustModel} \SelectedSample_{i} = \noNuclearPlants\\
    & \quad \sum_{i=1}^{\noSamplesinRobustModel} \sum_{t \in \TimeSteps}\outageSample{i}{t} \cdot \SelectedSample_{i} \leq \annualOutagePercentile, \label{eq:annual_unplanned_load_loss} \\
    & \quad \sum_{i=1}^{\noSamplesinRobustModel}\outageSample{i}{t} \cdot \SelectedSample_{i} \leq \maxNoOutagePercentile, \quad \forall t \in \AllHours \\
    & \quad \SelectedSample_{i} \in \left\{0, 1\right\}, \quad \forall i \in \left\{1, 2, \ldots, \noSamplesinRobustModel \right\}
\end{align}

\begin{figure*}[htbp]
    \centering
    \includegraphics[width=\linewidth]{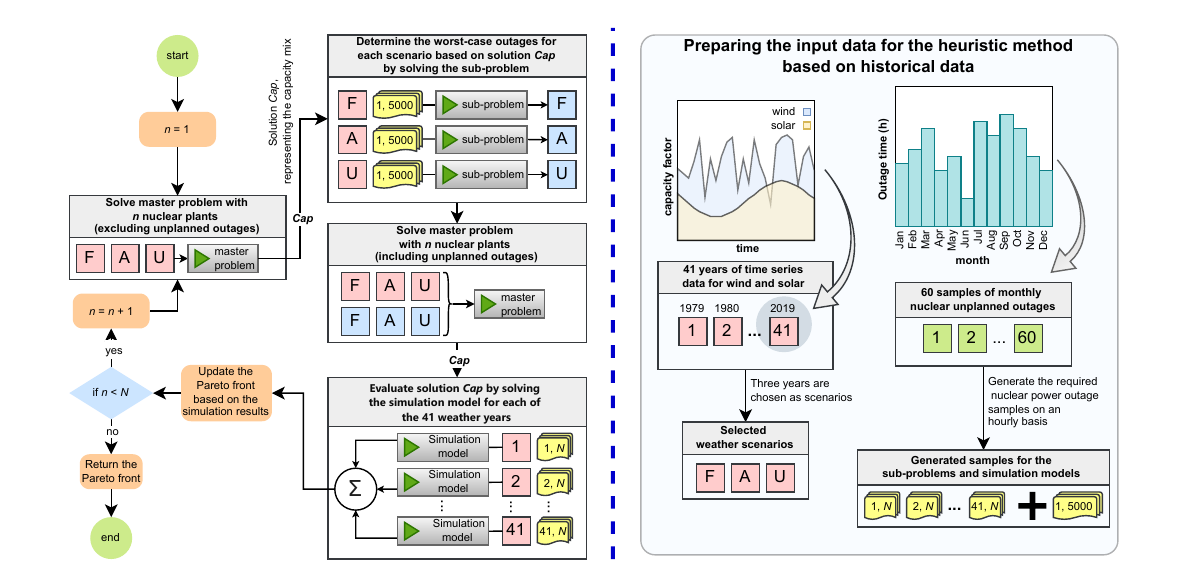}
    \caption{Flowchart of the heuristic method used to solve the \WNU (\Robust) model; \Favorable, \Average, and \Unfavorable represent \textit{Favorable}, \textit{Average}, and \textit{Unfavorable} weather scenarios, respectively.}
    \label{fig:flowchart}
\end{figure*}

The suggested heuristic method begins by solving the \hyperlink{Masterproblem}{Master problem}($\noNuclearPlants$, $\SelectedSampleBarBold$) for $\noNuclearPlants \in \Sizes$, where $\Sizes$ is the set of possible numbers of nuclear plants in the system. Since the optimal number of nuclear plants obtained from solving the \Robust model is unlikely to deviate significantly from that of the \Stochastic model, we can restrict $\Sizes$ by solving the \Stochastic model presented in Section~\ref{sec:StochasticModel}. To find a starting solution, for each $n \in \Sizes$, we set $\SelectedSampleBarBold = 0$ and solve the master problem to find a capacity mix, $\CapacityBarBold$. Then, for each scenario $s$, we calculate $\bar{\gamma}_{t}$ via Equation~\eqref{eq:normalized_net_load} and solve the \hyperlink{SubProblem}{Sub-problem} to obtain $\SelectedSampleBarBold$. Afterward, we solve the \hyperlink{Masterproblem}{Master problem} again to obtain an approximated solution for the \Robust model. Distributing the search over the set $\Sizes$ also allows us to incorporate simulation into the heuristic. Therefore, for the resulting capacity mix, we solve the \hyperlink{SimulationProblem}{Simulation problem} to evaluate the resulting solutions more accurately in terms of cost and loss of load. We can also examine different confidence levels, $\confidenceLevel$, to find a set of Pareto solutions. These processes are carried out for every $n$ in $\Sizes$ and $\confidenceLevel \in \mathcal{A}$, where $\mathcal{A}$ is the set of desired confidence levels. The flowchart of the heuristic method is presented in Figure~\ref{fig:flowchart}. 


The proposed method for solving the \Robust model, along with the associated mathematical models, was implemented using the Julia programming language. The pseudo-code of the heuristic, including its detailed steps, can also be found in \ref{sec:appendix_algorithms}. It should be noted that the optimization problems in this study are solved using the Gurobi\textregistered{} Optimizer \citep{gurobi}, with the barrier solver selected as the default.

\section{Experimental design} \label{sec:Experimentaldesign}

\begin{figure}[htbp]
    \centering
    \includegraphics[width=0.8\columnwidth]{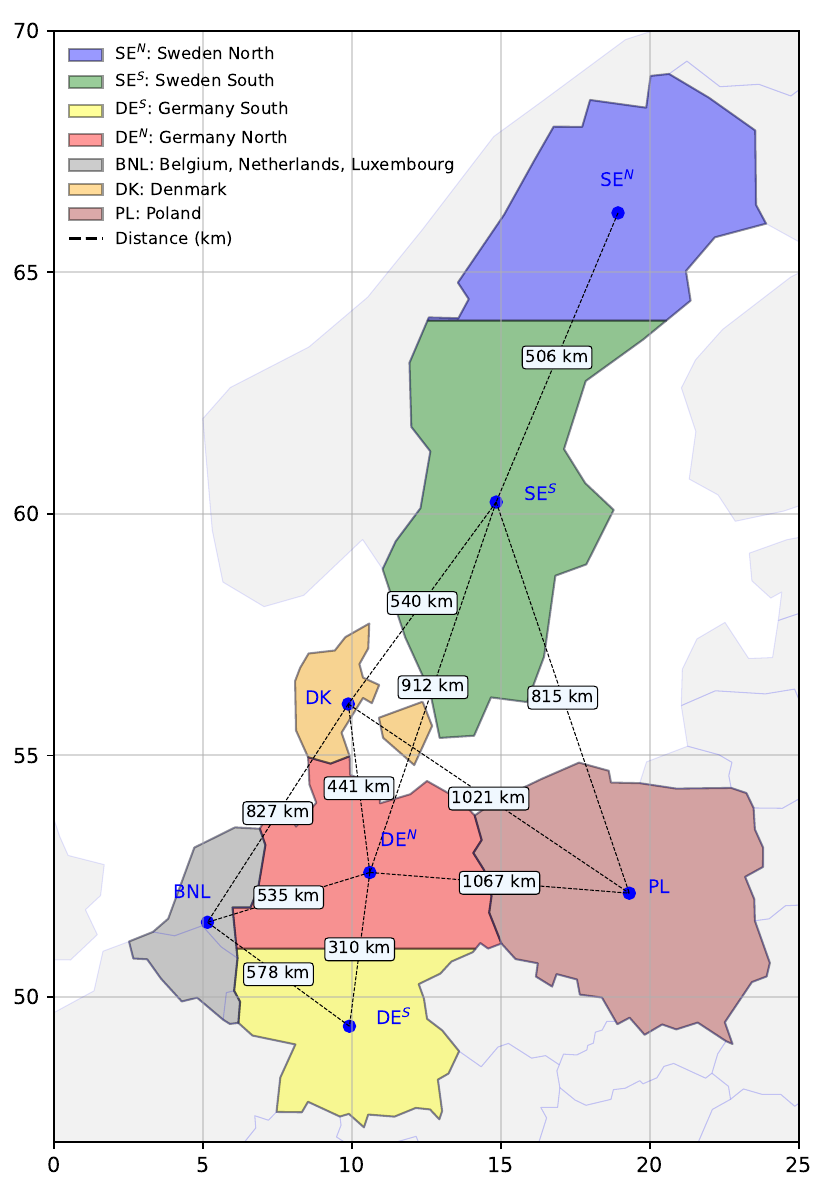}
    \caption{Selected regions in the case study, as well as the potential transmission lines among the regions.}
    \label{fig:regions_and_distances}
\end{figure}

The models presented in this study are tested in a case study for countries located in northern Europe. The case study includes seven regions: Northern Sweden (\SEN), Southern Sweden (\SES), Northern Germany (\DEN), Southern Germany (\DES), Poland (\PL), Denmark (\DK), and Belgium, the Netherlands, and Luxembourg (\BNL), grouped as one region. The data includes the time series of wind and solar capacity factors for 41 years, from 1979 to 2019. For the electricity load, projected time series data for the year 2050 is used. The potential links connecting the regions through transmission lines, along with the distances between the regions, $\distance{r}{r'}$, are shown in Figure~\ref{fig:regions_and_distances}. Through preliminary experiments, we found that a time step of seven hours, $\timeStep = \SI{7}{\hour}$, provides a suitable trade-off between accuracy and computational time. A description of the additional parameters in the case study can be found in \ref{sec:additional_parameters}.

\begin{figure*}[htbp]
    \centering
    \includegraphics[width=0.8\linewidth]{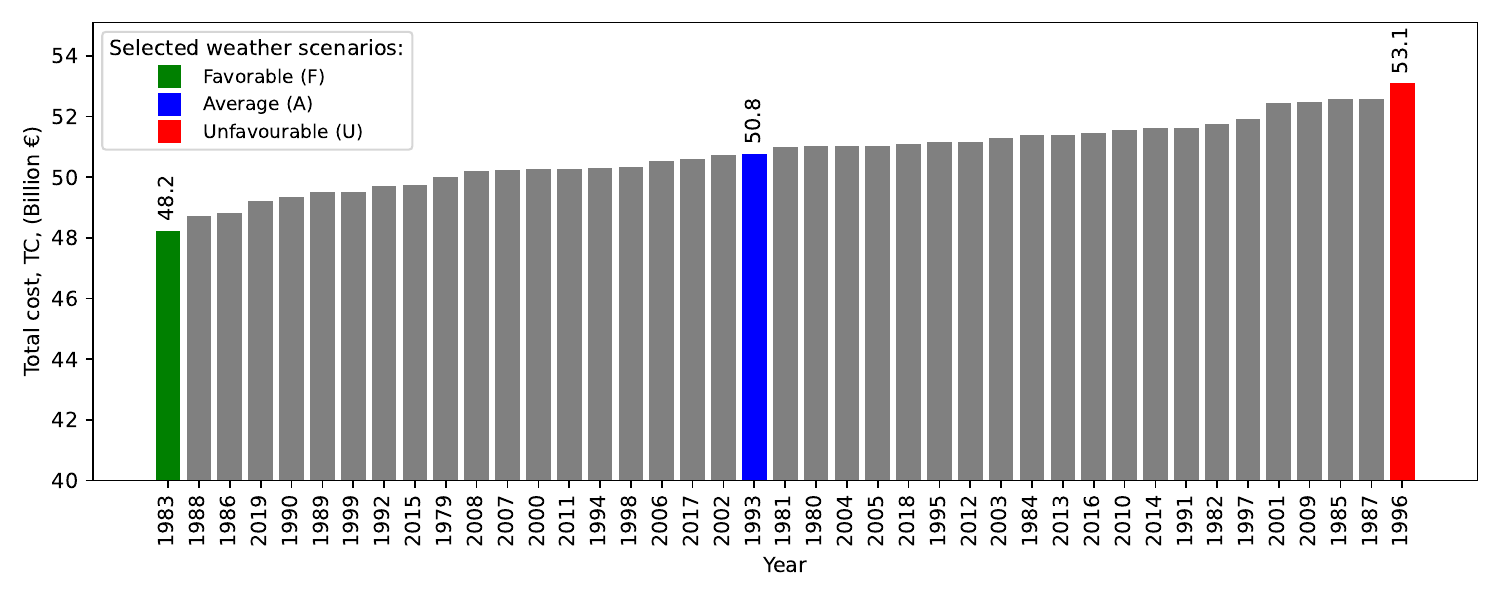}
    \caption{Optimal $\SC$ for each year, sorted in ascending order, with the three selected years highlighted to represent the weather scenarios.}
    \label{fig:LTC_individual_years}
\end{figure*}

To select weather scenarios, we conducted a preliminary test by solving 41 deterministic optimization problems (excluding unplanned nuclear power outages) based on data from 1979 to 2019. The objective function values from these tests, i.e., the $\SC$, serve as the basis for our scenario selection method. Figure~\ref{fig:LTC_individual_years} presents the years ranked according to their $\SC$ values. From this, three scenarios were chosen for the \WU and \WNU models: The \textit{Favorable} weather year scenario (\Favorable), based on the 1983 data (which results in the lowest \SC value); the \textit{Average} weather year scenario (\Average), based on the 1993 data (which has an \SC value closest to the average of the 41 $\SC$s); and the \textit{Unfavorable} weather year scenario (\Unfavorable), based on the 1996 data (which produces the highest $\SC$ value). It is assumed that scenario~\Favorable has a probability of 0.2, whereas both scenarios~\Average and \Unfavorable have probabilities of 0.4 each\footnote{These probabilities were determined through preliminary tests by solving the stochastic model and running simulations. The scenario combination that led to the lowest simulation cost was selected.}.

\begin{figure}[htbp]
    \centering
    \includegraphics[width=\columnwidth]{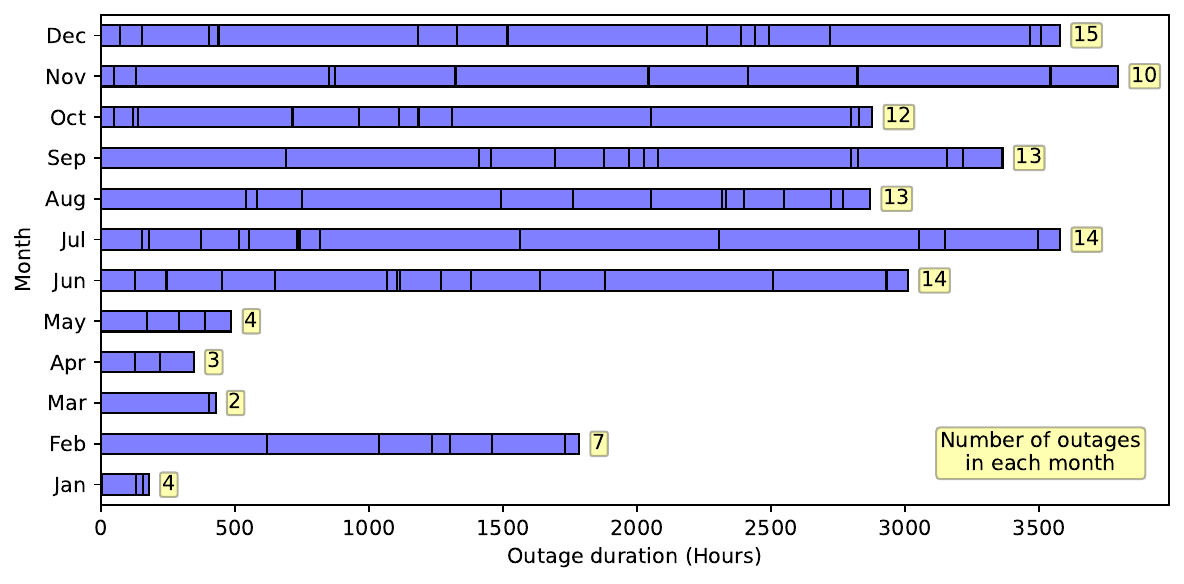}
    \caption{Unplanned nuclear power outages based on 60 separate 12-month samples from six nuclear plants in Sweden. Data source: \citep{iaea2023operating}.}
    \label{fig:Sweden_Unplanned_Outage}
\end{figure}

We have considered three models, presented in Table~\ref{tab:cases}, for comparison. For the \D (\Deterministic) model, only the normal scenario data is used (thus no weather uncertainty), and unplanned nuclear power outages are disregarded. In this case, we solve the \Stochastic model, presented in Section~\ref{sec:StochasticModel}, with only the given scenario~\Average. For the \WU (\Stochastic) model, the three scenarios, $\left\{\Favorable, \Average, \Unfavorable\right\}$, highlighted in Figure~\ref{fig:LTC_individual_years}, are used to account for uncertainty in weather conditions. Finally, the \WNU (\Robust) model is solved using the heuristic method presented in Section~\ref{sec:SolutionMethod} under three scenarios to account for uncertainties in unplanned nuclear power outages and weather conditions. In this respect, the heuristic method is employed with ten confidence levels, $\confidenceLevel \in \mathcal{A} = {0.1, 0.2, \ldots, 0.9} \cup {0.99}$, and 13 possible sizes for the number of nuclear plants in the system, $\noNuclearPlants \in \Sizes = \left\{36, 38, \ldots, 62\right\}$. In the experiments, each nuclear plant is assumed to have a capacity of 1~\unit{\GWh} ($\oneNuclearPlantSize = 1000$). Additionally, we set $\noSamplesinRobustModel = 50000$ and $\noSamplesforPercentiles = 100000$ (see Sections~\ref{sec:RobustModel} and \ref{sec:SolutionMethod}, for further details about these two parameters). Historical monthly outage data from six nuclear power plants in Sweden are used in this study (see Figure~\ref{fig:Sweden_Unplanned_Outage}).

The capacity layouts resulting from the models are evaluated using simulations. The following are the four simulation cases used for this purpose:

\begin{itemize} 
    \item \hypertarget{No-outage}{\NoOutageSimulation} simulation: 41 historical weather years, neglecting the risk of unplanned nuclear power outages. 
    \item \hypertarget{Normal}{\NormalSimulation} simulation: 41 historical weather years, including the risk of unplanned nuclear power outages. 
    \item \hypertarget{Unfavorable weather}{\UnfavorableSimulation} simulation: 41 extreme weather years (with historical data from the year 1996 used for all 41 weather years), combined with the risk of nuclear power outages. This simulation increases the likelihood of nuclear unavailability coinciding with unfavorable weather and demand situations. 
    \item \hypertarget{Dunkelflaute}{\DunkelflauteSimulation} simulation: Similar to the normal simulation, but with an additional \emph{Dunkelflaute} event starting on February 1st (day 32 of the 365-day year) and lasting for 14 days. \footnote{This was done by modifying the wind and solar capacity factors within a specific period. During this period, the original wind and solar capacity factors of the regions in each simulation year were multiplied by a constant, referred to as Dunkelflaute intensity, which was assumed to be 0.4.}. \end{itemize}

The solutions of the \Deterministic and \Stochastic models are evaluated using the \NoOutageSimulation and \NormalSimulation simulation cases (see Section~\ref{sec:stochastic_results}), while the last three simulation cases are used to evaluate the solutions of the \Robust model (see Section~\ref{sec:robust_results}).

\begin{table}[htbp]
    \centering
    \caption{Models considered for making comparisons in this paper.}
    \label{tab:cases}
    \resizebox{\columnwidth}{!}{
    \begin{threeparttable}   
        \begin{tabular}{llcccc} \toprule
        Model & Description & \begin{tabular}[c]{@{}c@{}}Unplanned nuclear \\ power outages\end{tabular} & \multicolumn{3}{l}{Scenarios\tnote{a}} \\ \midrule
        \Deterministic\tnote{b} & Deterministic (no uncertainty) &  &  & \Average &  \\
        \Stochastic\tnote{c} & Uncertainty in weather conditions &  & \Favorable & \Average & \Unfavorable \\
        \Robust\tnote{c} & \begin{tabular}[l]{@{}l@{}}Uncertainty in weather conditions \\ and nuclear power availability\end{tabular}   & \checkmark & \Favorable & \Average & \Unfavorable \\ \bottomrule
        \end{tabular}
    \begin{tablenotes}
        \item[a] Scenario are as follows: \Favorable: Favorable weather year, \Average: Average weather year, and \Unfavorable: Unfavorable weather year. 
        \item[b] In this case, the probability of Scenario~\Average is 1.   
        \item[b] In this case, the probabilities of the scenarios are as follows: $\prob{\Favorable} = 0.2$ and $\prob{\Average} = \prob{\Unfavorable} = 0.4$.
    \end{tablenotes}
    \end{threeparttable}}
\end{table}
\subsection{Results}\label{sec:Results}

In Section~\ref{sec:stochastic_results}, we begin by comparing the solutions obtained from the \Deterministic and \Stochastic models, which are similar to models already published in the literature, in order to investigate the effect of uncertainty in weather conditions and evaluate how the system performs when unplanned nuclear power outages are disregarded. Then, in Section~\ref{sec:robust_results}, we proceed to investigate the performance of our proposed \Robust model, which takes into account uncertainty regarding nuclear availability in the planning phase.

The energy systems resulting from the optimization in the \Deterministic and \Stochastic models are depicted in \ref{sec:appendix_figures}, Figure~\ref{fig:D_and_WU_map}. The system is primarily characterized by generation from wind and solar power, with nuclear energy closely following. The transmission capacity is significantly higher than current levels, and gas serves as the primary backup capacity. Additionally, the battery capacity is designed to cover approximately two hours of average demand.

\subsubsection{The Effect of Considering Weather Uncertainty: Comparing the Performance of the \Deterministic and \Stochastic Models}\label{sec:stochastic_results}

The main purpose of this section is to demonstrate that failing to consider unplanned nuclear power failures during the planning stage is consequential if the goal is to design a robust energy system with minimal loss of load. To this end, we begin by examining the performance (total cost and Loss of Load (LoL)) of the \Deterministic model, which has no uncertainty, and the \Stochastic model, which incorporates uncertainty in weather conditions. We achieve this by exposing the respective solutions to the \NoOutageSimulation and \NormalSimulation simulation environments. The difference in system cost between the performance of the solution from the \Stochastic model and that from the \Deterministic model when exposed to these simulation environments is termed the \emph{value of stochastic solution} (VSS) \citep{birge2011introduction, backe2021stable}. The simulation results indicate that exposing the solutions only to varying weather conditions (\NoOutageSimulation simulation) yields a VSS of $(51309 - 51254) = 55$~\unit{\MCurrency}, which is equivalent to annual savings of 0.11\%.

\begin{table*}[htbp]
\centering
    \caption{Summary of the computational results for the \D (\Deterministic) and \WU (\Stochastic) models.} \label{tab:result_summary}
    \resizebox{\textwidth}{!}{
      \begin{tabular}{@{}lccccccccc@{}}
            \toprule
                                         & \multicolumn{2}{c}{\begin{tabular}[c]{@{}c@{}}Results of solving \\ the models without simulation \end{tabular} } 
                                         & \multicolumn{2}{c}{\begin{tabular}[c]{@{}c@{}}Simulation results excluding \\ unplanned nuclear power outages \\ (\NoOutageSimulation simulation) \end{tabular}} 
                                         & \multicolumn{2}{c}{\begin{tabular}[c]{@{}c@{}}Simulation results with \\ unplanned nuclear power outages \\ (\NormalSimulation simulation) \end{tabular}} \\ \cmidrule(lr){2-3} \cmidrule(lr){4-5} \cmidrule(lr){6-7}
                                         & \Deterministic         & \Stochastic           
                                         & \Deterministic         & \Stochastic
                                         & \Deterministic         & \Stochastic\\ \midrule
            $\SC$ (\unit{\MCurrency})                     & 50777   & 51534   & 51309   & 51254   & 52301   & 52205   \\ 
            LoL (\unit{\GWh})                           & 0       & 0       & 13      & 7       & 34      & 18      \\
            LoL frequency                               & $0/41$  & $0/41$  & $7/41$  & $7/41$  & $14/41$ & $14/41$   \\
            LoL percentage (\%)                         & 0       & 0       & $1.3 \times 10^{-3}$  & $6.6 \times 10^{-4}$     & $3.3 \times 10^{-3}$   & $1.7 \times 10^{-3}$   \\
            Investment cost (\unit{\MCurrency})            & 29539   & 30245   & 29539   & 30245   & 29539   & 30245   \\ 
            Fixed cost (\unit{\MCurrency})                 & 12064   & 12292   & 12064   & 12292   & 12064   & 12292   \\ 
            Operational cost (\unit{\MCurrency})           & 6262    & 6331    & 6636    & 6239    & 7092    & 6689    \\ 
            Emissions cost (\unit{\MCurrency})             & 2702    & 2514    & 2978    & 2440    & 3449    & 2905    \\ 
            Load shedding cost (\unit{\MCurrency})         & 210     & 152     & 92      & 38      & 157     & 74      \\ \bottomrule
    \end{tabular}     }
\end{table*}

Regarding Loss-of-Load (LoL), the solution of the \Deterministic model exhibits a LoL that is twice as high as the solution of the \Stochastic model in the \NoOutageSimulation simulation. When there is uncertainty in both nuclear power outages and weather conditions (\NormalSimulation simulation), the value of stochastic solution (VSS) is calculated as $52301 - 52205 = 96$\unit{\MCurrency}, or $0.18$\%. In this case, the LoL remains approximately twice as high in the \Deterministic model’s solution compared to the \Stochastic model’s solution ($34$ vs. $18$ \unit{\GWh}). These results demonstrate that, although the \Stochastic model performs better than the \Deterministic model, the solution is not sufficiently robust, as LoL increases significantly when nuclear power outages are considered. As seen in Table~\ref{tab:result_summary}, LoL also becomes more frequent with unplanned nuclear power outages, rising from 7 to 14 events over the 41 years. Detailed simulation results for the solutions under the \NormalSimulation simulation case can be found in Table~\ref{tab:detailed_simulation_results}.

To assess how the vulnerability of the energy system varies with nuclear capacity, we conducted a sensitivity analysis on the number of nuclear plants. The results revealed that as nuclear power penetration increases, the system becomes more susceptible to loss of load (LoL), as shown in Figure~\ref{fig:LoL_plot_different_n} in \ref{sec:appendix_figures}. Furthermore, the Pareto front analysis (Figure~\ref{fig:pareto_front_stochastich_and_deterministic} in \ref{sec:appendix_figures}) indicates that the \Stochastic model outperforms the \Deterministic model in terms of system cost and LoL. Based on these findings, the following section will use the Pareto front derived from the \Stochastic model to evaluate the solutions of the \Robust model.
\subsubsection{Considering nuclear uncertainty in the planning phase: Simulation results for the \Robust) model} \label{sec:robust_results}

\begin{figure*}[htbp]
    \centering
    \includegraphics[width=\linewidth]{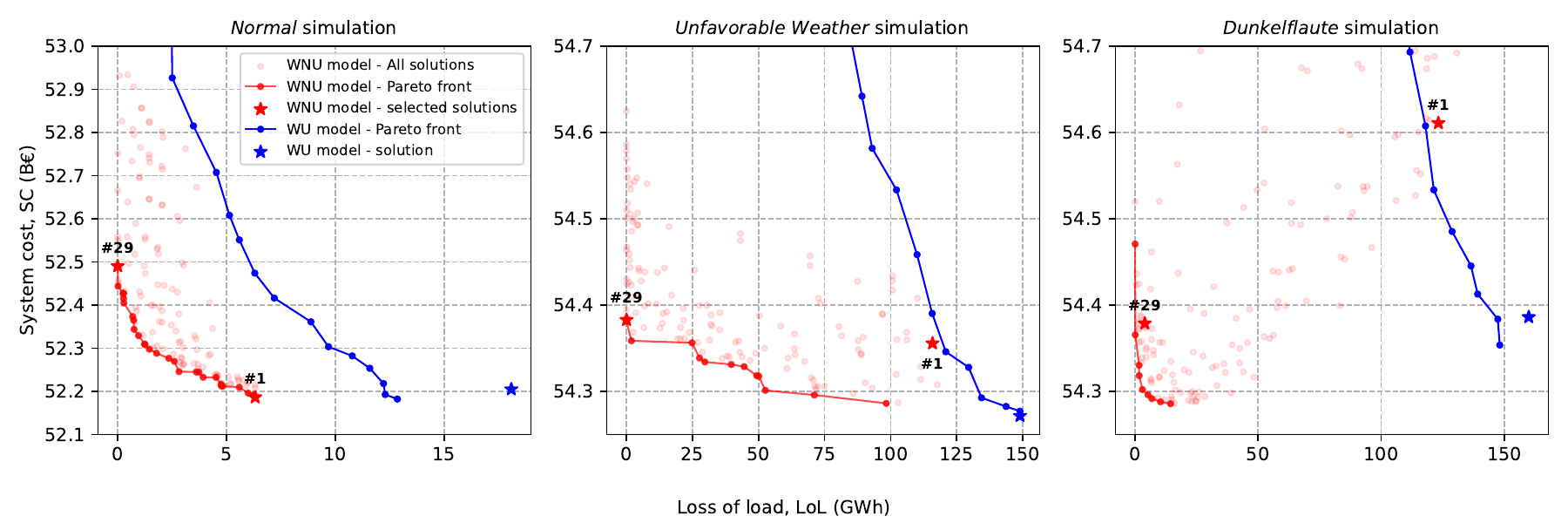}
    \caption{Pareto front comparing the solutions of robust and stochastic models across three simulation cases; for the \Stochastic model, the Pareto front has been obtained through a sensitivity analysis on the number of nuclear plants in the systems.}
    \label{fig:Pareto_front_comparison_all}
\end{figure*}

The \WNU (\Robust) model incorporates uncertainty in both weather conditions and nuclear power outages during the planning phase. Its performance was evaluated under three simulation conditions: \NormalSimulation, \UnfavorableSimulation, and \DunkelflauteSimulation, with the latter two scenarios increasing the likelihood of unfavorable weather conditions coinciding with unplanned nuclear power outages. The corresponding Pareto fronts for the \Stochastic and \Robust models under each simulation case are shown in Figure~\ref{fig:Pareto_front_comparison_all}\footnote{In this figure, the Pareto front for the \Stochastic model is derived from sensitivity analysis on the number of nuclear power plants, as discussed in \ref{sec:sensitivity analysis} and illustrated in Figure~\ref{fig:LoL_plot_different_n}.}. The figure demonstrates that the \Robust model consistently outperforms the \Stochastic model in both loss-of-load (LoL) and system cost (\SC) across all simulation cases. Notably, in all three simulation conditions, the \Robust model produces solutions with zero LoL, underscoring its suitability for designing a more robust energy system.

Focusing on the \NormalSimulation simulation case in Figure~\ref{fig:Pareto_front_comparison_all}, the results for the \Robust model reveal a distinct Pareto front with 29 solutions, showing increasing cost to achieve zero LoL. The total capacity mix, costs, and LoL for all 29 solutions are detailed in Table~\ref{tab:pareto_solutions_summary} in the appendix. The two endpoints of the Pareto front, representing the most and least conservative solutions, are highlighted with stars in Figure~\ref{fig:Pareto_front_comparison_all}\footnote{For these solutions, the simulation results for each individual year, along with the capacity layout, are available in Table~\ref{tab:detailed_simulation_results} and Figure~\ref{fig:robust_solutions}.}.

The simulation results indicate a substantial difference in LoL across different years, see \ref{tab:detailed_simulation_results} in \ref{sec:appendix_tables}, but here we show the average results. Based on the two solutions obtained by the \Robust model (solutions \#1 and \#29), the additional cost that we may be willing to pay to further protect the energy system against LoL, the \emph{price of robustness} \citep{bertsimas2004price}, is calculated as $52490 - 52186 = 304$~\unit{\MCurrency}, which is equivalent to a 0.6\% increase in the system cost.

In the simulation cases \UnfavorableSimulation and \DunkelflauteSimulation shown in Figure~\ref{fig:Pareto_front_comparison_all}, the solutions of the \Robust and \Stochastic models are exposed to a more harsh environment, both because the weather situation is less favorable, but also because there is a greater risk of unfavorable weather and demand coinciding with nuclear plants being unavailable. These conditions give rise to a higher system cost, by about 4\%, and a higher value for LoL, than in the \NormalSimulation case, see Figure~\ref{fig:Pareto_front_comparison_all}. Yet, there is one solution, namely the solution that produced zero LoL in the simulation case \NormalSimulation (solution \#29), that yields zero LoL also in the simulation cases \UnfavorableSimulation and \DunkelflauteSimulation. We may thus infer that the \Robust model is capable of finding solutions that are resilient not only to the conditions under which it was optimized but also, more generally, to unfavorable combinations of weather, demand, and nuclear power outages.

\begin{figure}[htbp]
    \centering
    \includegraphics[width=\columnwidth]{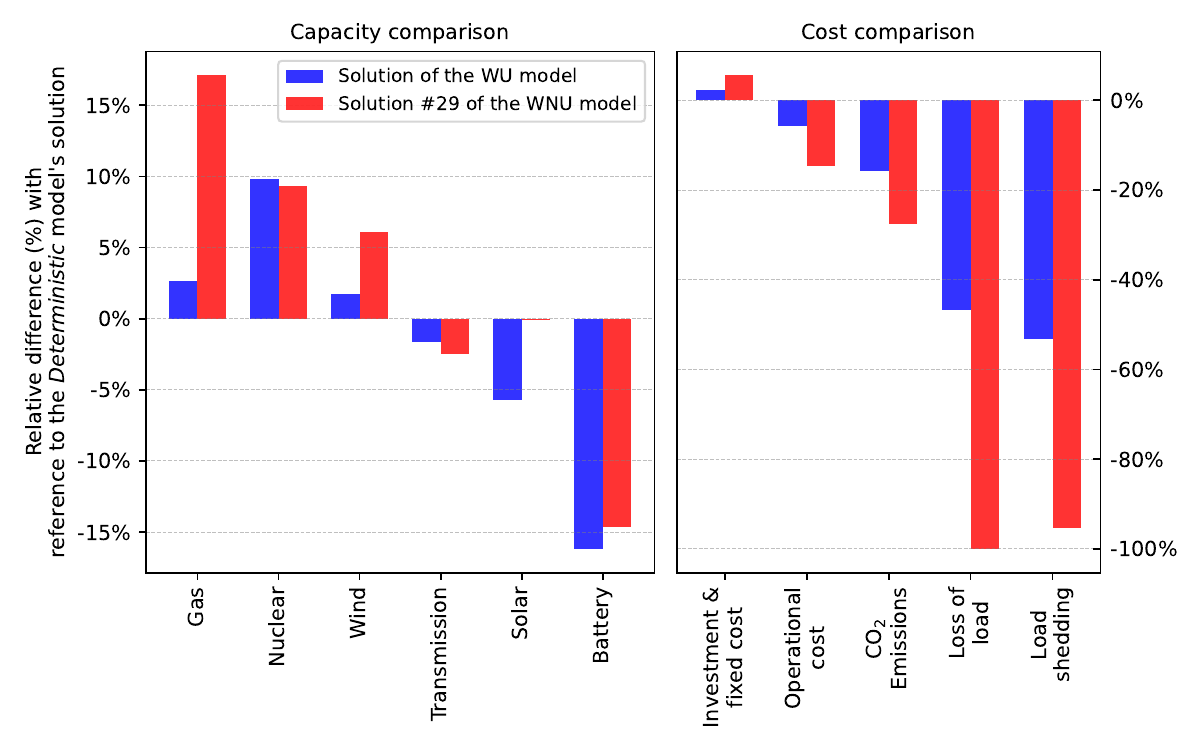}
    \caption{Relative differences in capacity and costs, comparing the solutions of the \WU (\Stochastic) and \WNU (\Robust) models against the solution of the \D (\Deterministic) model; the cost comparison is based on the \NormalSimulation simulation case.}
    \label{fig:relative_comparion_capacity_and_cost}
\end{figure}

Figure~\ref{fig:relative_comparion_capacity_and_cost} shows the change in capacity mix (LHS) and performance measures (RHS) between the \Stochastic solution, the most conservative solution of the \Robust model (solution \#29), and the \Deterministic solution. As seen in the figure, compared to the solution of the \Deterministic model, the \Robust solution is mainly characterized by i) an increase in gas, nuclear, and wind capacity (by 17\%, 9\%, and 6\%, respectively), and ii) a decrease in battery capacity (by 15\%). Figure~\ref{fig:relative_comparion_capacity_and_cost} also shows that both LoL and load shedding decrease compared to the \Deterministic model. Thus, the change in capacity mix dominated by added gas capacity provides a greater ability to avoid both load-shedding and LoL. Yet, there is less total generation by gas, which is also evident from the lower (28\%) \ch{CO2} emissions in the \Robust model's solution. Note that nuclear capacity increases by around 10\% in the solutions of the \Stochastic and \Robust models compared to the \Deterministic model's solution, see Figure~\ref{fig:relative_comparion_capacity_and_cost}. In other words, robustness is provided not only with gas, but also with added nuclear capacity, even when, as is the case in the \Robust formulation, nuclear has the disadvantage of unplanned outages.

\begin{figure}[htbp]
    \centering
    \includegraphics[width=0.8\columnwidth]{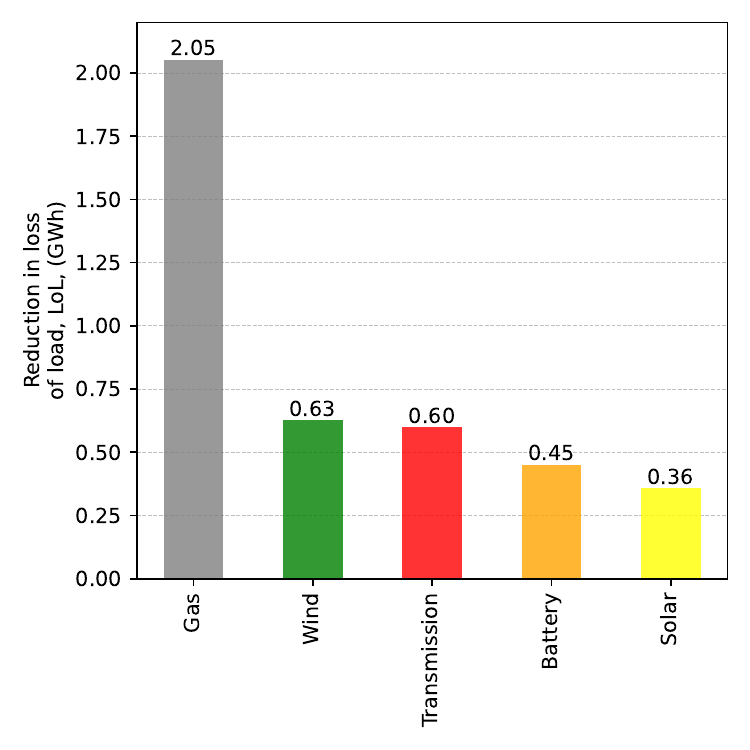}
    \caption{Impact of a one-unit change in capacities on the reduction of LoL.}
    \label{fig:LR_coefficient}
\end{figure}

For each $\noNuclearPlants \in \Sizes$, we can subtract the capacities and LoL in the \Robust model's solutions from those in the corresponding \Stochastic model's solutions. Since each solution pair has the same nuclear power capacity, this adjustment isolates the impact of other technologies. We then fit a linear regression model to the adjusted data to evaluate how other power generation sources contribute to LoL reduction\footnote{The regression analysis yielded a multiple $R$ of 0.983 ($R^2 = 0.966$), indicating a strong positive correlation between the adjusted capacity values (predictors) and the adjusted LoL (response variable). This suggests that 96.6\% of the variance in LoL reduction is explained by the capacity differences of the non-nuclear technologies. The statistical significance of the regression model is further supported by a near-zero $p$-value (the $F$-statistic of the regression model was about 770).}. The regression coefficients for each technology are plotted in Figure~\ref{fig:LR_coefficient}, in decreasing order. This result demonstrates approximately how much an increase in the capacity of each technology contributes to LoL reduction. As expected, gas technology, with a regression coefficient of almost two, plays a major role in improving the robustness of energy systems; in other words, an additional 1~\unit{\GW} of gas capacity reduces LoL by 2~\unit{\GWh}. Following gas, wind, transmission, battery, and solar technologies contribute to LoL reduction in decreasing order of impact. 

\begin{figure}[htbp]
    \centering
    \includegraphics[width=\columnwidth]{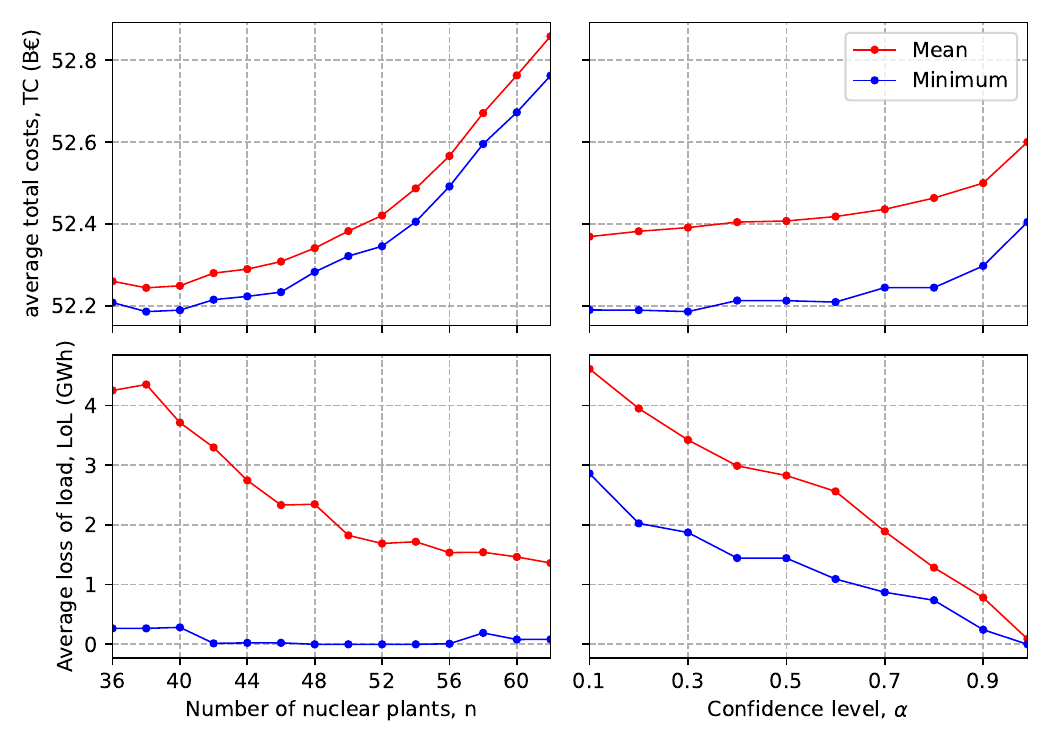}
    \caption{Effect of confidence level and number of nuclear plants on the \WNU (\Robust) model's solutions.}
    \label{fig:effect_of_alpha_and_n_2D}
\end{figure}

As explained earlier, the \Robust model was solved for ten different confidence levels ($\confidenceLevel$ ranging from 0.1 to 0.99) and 13 possible sizes for the number of nuclear plants in the system ($\noNuclearPlants$ ranging from 36 to 62 nuclear plants), the results consisting of 29 Pareto solutions, are provided in Table~\ref{tab:pareto_solutions_summary}). Figure~\ref{fig:effect_of_alpha_and_n_2D} shows the effect of the two parameters (confidence level and number of nuclear plants) on the cost and LoL. In this figure, the red lines represent the mean values and the blue lines represent the minimum values across the number of nuclear plants (left panels in Figure~\ref{fig:effect_of_alpha_and_n_2D}) and the confidence level (right panels in Figure~\ref{fig:effect_of_alpha_and_n_2D}). The result in the lower left panel in Figure~\ref{fig:effect_of_alpha_and_n_2D} indicates that the impact of the number of nuclear plants ($\noNuclearPlants$) on the minimum LoL is minimal. This suggests that regardless of the number of nuclear plants, the \Robust model consistently identifies solutions that are resilient to LoL. On the other hand, as shown in the top left panel of Figure~\ref{fig:effect_of_alpha_and_n_2D}, both the mean and minimum $\SC$ increase as $\noNuclearPlants$ grows; this is due to the higher costs associated with expanding nuclear power capacity. In the \Robust model, the confidence level, $\alpha$, is the main parameter that adjusts the conservatism of the \Robust model. As expected, increasing $\alpha$ leads to a higher $\SC$ and a lower LoL for both the mean and minimum values, see the two right panels in Figure~\ref{fig:effect_of_alpha_and_n_2D}.

\subsection{Discussion}\label{sec:Discussion}

This paper presents the first proof of concept for a novel model, named \WNU (\Robust) model, to account for uncertainties in both weather conditions and nuclear power availability in an electricity system. While numerous studies have investigated the importance of weather uncertainty, either through stochastic programming or by comparing solutions based on different weather years \citep{ringkjob2020short, seljom2015short, gotske2024designing, collins2018impacts}, the uncertainty of nuclear power output has not been well investigated in the energy system literature. We show that failing to account for the uncertainty arising from nuclear power plant outages is consequential.
\subsubsection{Method contribution and comparison}

The proposed \Robust model utilizes a data-driven robust optimization approach to identify optimal investment decisions under the uncertainties as mentioned earlier. The results of our model were evaluated using two metrics: system cost and Loss-of-Load (LoL). We introduced LoL as a measure to represent situations where the system cannot meet demand through load shedding or other standard measures, resulting in a risk of brownouts or blackouts. The choice to show LoL as a separate quantity, and not just as an added cost, also underscores the trade-off between, on the one hand, system cost and, on the other, system resilience or robustness. As our model does not detail the power system, we are unable to analyze the precise consequences or the economic impact of instances of load deficiency. Instead, LoL is used as an indicator to reflect the risk of high-cost events caused by electricity supply shortages. The model can lead to solutions that are robust in relation to power shortages while at the same time being cost-effective, as shown in Figure~\ref{fig:Pareto_front_comparison_all}. 

The designated heuristic to solve the \WNU model is computationally efficient enough to be incorporated into an energy system investment model with seven nodes, approximately 1,250 time steps, and eight generation technologies, in addition to storage and transmission. In terms of performance, it has the ability to decrease LoL to zero or near zero even for a high penetration level of nuclear power, while at the same time minimizing total costs. 
\subsubsection{The value and cost of robust design, and comparison to literature}

Energy system models, including those used in papers that focus on the role of nuclear power \citep{hjelmeland2025role, duan2022stylized, thellufsen2024cost, price2023role, fattahi2022analyzing, sepulveda2018role, kan2020cost}, represent it as a stable, dispatchable power source, without accounting for the risk of unplanned outages. Here, we show that representing uncertainty in nuclear power availability is important in the investment phase to increase energy system resilience by avoiding LoL: When nuclear uncertainty was considered, i.e. when comparing the results from the solution of the \WNU model with those of the solution of the \WU model, LoL decreased from 6~\unit{\GWh} to zero under normal conditions (\NormalSimulation simulation), see Figure~\ref{fig:Pareto_front_comparison_all}. Furthermore, when the weather conditions were worse (\UnfavorableSimulation and \DunkelflauteSimulation simulation), our proposed model still managed to find solutions with near-zero LoL, whereas the \WU model was unable to achieve this cost-effectively.

Our results, as well as those in other papers, such as \citet{seljom2015short,seljom2021stochastic}, thus highlight the value (in terms of lower cost, lower LoL or both) of designing systems that are robust to varying weather conditions and, in our case, nuclear power outage conditions. However, building for robustness comes with an initially higher cost, because the investment part of the system cost is higher than that estimated by the corresponding deterministic model. We find that this additional projected cost, i.e., the total system cost projected by the investment model, that is resilient under both weather and nuclear uncertainty (solution produced by the \WNU model), is 6\% higher compared to the estimate of the \D model. The additional cost to build a system that considers only weather uncertainty (solution from \WU model) is 2\% higher in this study, while previous studies reported a 1\% higher cost \citep{seljom2021stochastic} and up to 20\% higher \citep{ringkjob2020short}.  

Our results also indicate that the main impact on the capacity mix from incorporating weather and nuclear uncertainties is an increased capacity of gas turbines. A previous paper that investigated weather variations \citep{gotske2024designing} \emph{assumed} that adding gas turbines was the technology that would provide increased resilience to the energy system, but here we show that this is indeed the case. This result is also consistent with the findings of \citet{seljom2015short, seljom2021stochastic}, which demonstrate that accounting for weather uncertainty leads to an increase in the dispatchable thermal generation capacity. Interestingly, the effect on nuclear investments from accounting for nuclear uncertainty is not, as one might expect, a reduction in nuclear power investments but rather a small increase, see Figure \ref{fig:relative_comparion_capacity_and_cost}. Further investigations are needed in order to determine whether this result is stable also under varying surrounding assumptions.
\subsubsection{Limitations: data and generalization to larger models}

The exact effect of nuclear outage risks is challenging to estimate at this stage, as we used Swedish data as a proxy. Whether this data accurately reflects the performance of \emph{future} nuclear power plants remains to be assessed: it may be the case that the reasons behind the outage statistics are associated with an aging power plant fleet. Moreover, the impact of nuclear power availability likely varies depending on the overall energy system configuration. Nevertheless, our findings highlight that unplanned outages should be considered when evaluating the resilience of energy systems with nuclear power — a perspective that is largely absent in the existing literature. 

The test case for our model is an overnight investment capacity expansion model with seven nodes, 1250 time steps, and eight technology options. The size and design of this model are simpler than the current state-of-the-art models, such as PyPSA \citep{horsch2018pypsa} or GenX \citep{jenkins2017enhanced}, which, apart from being larger in size, also incorporates more technology options, and includes sectors with demands for hydrogen. We believe that our formulation would be feasible to incorporate also in these larger models, but we will leave it for future work to investigate the effect of nuclear power outages in such settings.
\subsubsection{Implications for modeling and future work}

Our results indicate that the solutions generated by the proposed \Robust model, which accounts for both nuclear and weather uncertainties, outperform those that consider only weather uncertainty, particularly during extreme weather events such as Dunkelflaute. This finding indicates that the robust measures implemented to address unplanned nuclear power outages also enhance the system's ability to handle other stresses effectively. The future energy system faces numerous uncertainties (generation, transmission, demand, etc.), and representing and modeling all of them is a challenging task, as stochastic models can become computationally intensive. Our results show that addressing two key uncertainties enhances the robustness of the system. However, identifying other important uncertainties to further increase the system's robustness remains a crucial area for future research.

In this study, we focus mainly on uncertainty within \emph{normal} conditions and do not account for extreme events such as storms that damage infrastructure, terrorist attacks, war, or similar risks. Addressing such extreme events may require a different set of measures to achieve robustness. Designing a robust energy system that accounts for these risks, as well as identifying the parameters that are most crucial for added robustness, remains an important area for future research.
\section{Conclusion}\label{sec:Conclusion}

This paper investigated a capacity expansion problem for electricity with uncertain weather conditions and nuclear power availability. A scenario-based stochastic optimization approach was employed to represent varying weather conditions. The stochastic model was then extended to incorporate unplanned nuclear power outages using a data-driven adjustable robust optimization approach. Given the computational complexity of the robust problem, a heuristic method was proposed to efficiently solve the model. Since the stochastic and robust models were based on a limited set of scenarios chosen from historical weather data, simulations were carried out using a larger set of historical data to evaluate the performance of the models more accurately under different conditions. We tested the models in a case study for Northern Europe, consisting of seven nodes, with a time resolution of 1250 time steps. Our findings show that:
\begin{itemize}
    \item The current standard energy system model formulations, both deterministic and those that incorporate uncertainty in weather, are vulnerable to Loss-of-Load (LoL) when exposed to nuclear power uncertainty. This effect is exacerbated at high penetration levels of nuclear power. 
    \item The robust model that incorporates uncertainty in nuclear power availability yielded solutions with no or very low levels of loss of load (LoL), even at high levels of nuclear penetration. These solutions also performed well when tested against new and more extreme combinations of weather events (e.g., Dunkelflaute) and nuclear power outages.
    \item The price of robustness (the cost to build a system that suffers from less LoL) in our test case was 304~\unit{\Currency \per \MWh}, which corresponds to a 0.6\% increase in the system costs.
    \item Investment strategies that provide added robustness are dominated by greater gas capacity. In addition, and perhaps counterintuitive, the introduction of variability and uncertainty of nuclear operation did not significantly change the optimal amount of nuclear power. Rather, robust solutions showed a slightly higher optimal capacity of nuclear power, compared to the deterministic solution.
\end{itemize}

Based on the findings of this paper, we conclude that the uncertainty in nuclear power availability in European countries may justify its inclusion in energy system models that inform policies on energy system resilience. Future research should aim to offer more concrete guidance on resilient strategies for countries seeking to incorporate nuclear power into their \ch{CO2}-neutral energy mix. 
\section*{Acknowledgments}

We would like to thank Iegor Riepin for valuable discussions during the initial stage of this project. This work was funded by the Swedish Energy Agency under grant numbers P2023--01323 and P2022--00768.

\bibliographystyle{elsarticle-num-names} 
\bibliography{cas-refs}
\appendix
\setcounter{figure}{0}
\setcounter{table}{0}
\section{Notation}\label{sec:appendix_notation}

\nomenclature[A01]{$\AllHours$}{Set of hours in a year ($\AllHours = \{ 1, 2, \ldots, 8760 \}$)}
\nomenclature[A02]{$\TimeSteps$}{Set of all time steps ($\TimeSteps = \{ t \in \AllHours \mid h = 1 + k \cdot \timeStep, \, k \in \mathbb{N}_0 \}$)}
\nomenclature[A03]{$\Scenarios$}{Set of scenarios used in the robust and stochastic models}
\nomenclature[A04]{$\ScenariosSimulation$}{Set of scenarios used for simulations}
\nomenclature[A05]{$\Regions$}{Set of regions}
\nomenclature[A06]{$\Plants$}{Set of power generation plants ($\Plants = \{\Gas, \PV, \Wind, \Hydro, \Nuclear\}$)}
\nomenclature[A07]{$\Technologies$}{Set of technologies ($\Technologies = \Plants \cup \{\Battery, \Inverter\}$)}
\nomenclature[A08]{$\ValidLinks$}{Set of valid transmission power lines that connect regions}

\nomenclature[B01]{$s \in \Scenarios$}{Index for scenarios}
\nomenclature[B02]{$p \in \Technologies$}{Index for technologies}
\nomenclature[B03]{$t \in \TimeSteps$}{Index for hours and time steps}
\nomenclature[B04]{$r, r' \in \Regions$}{Index for regions}

\nomenclature[C01]{$\timeStep$}{Duration of each time period \nomunit{\unit{\hour}}}
\nomenclature[C02]{$\cf{s}{r}{p}{t}$}{Capacity factor of plant $p$ in region $r$ under scenario $s$ at time step $t$}
\nomenclature[C03]{$\plannedOutageRate$}{Planned outage rate for nuclear plants}
\nomenclature[C04]{$\hydroReservoir{r}$}{Maximum hydro reservoir level in region $r$ \nomunit{\unit{\MWh}}}
\nomenclature[C05]{$\resInflow{r}{t}$}{Water inflow to reservoir in region $r$ at time step $t$ \nomunit{\unit{\MWh}}}
\nomenclature[C06]{$\distance{r}{r'}$}{Distance between regions $r$ and $r'$ \nomunit{\unit{\km}}}
\nomenclature[C07]{$\load{r}{t}$}{Electricity demand in region $r$ at time step $t$ \nomunit{\unit{\MW}}}
\nomenclature[C08]{$\transmissionEfficiency{r}{r'}$}{Transmission efficiency of electricity between $r$ and $r'$}
\nomenclature[C09]{$\sheddingCapRatio$}{Maximum fraction of electricity demand that can be shed per each time step}
\nomenclature[C10]{$\emissionFactor{p}$}{Emission factor of plant $p$ \nomunit{\unit{\TonCOtwo\per\MWh}}}
\nomenclature[C11]{$\efficiency{p}$}{Efficiency of technology $p$}
\nomenclature[C12]{$\fuelCost{p}$}{Fuel cost of technology $p$ \nomunit{\unit{\Currency\per\MWh}}}
\nomenclature[C13]{$\varCost{p}$}{Variable cost of technology $p$ \nomunit{\unit{\Currency\per\MWh}}}
\nomenclature[C14]{$\fixedCost{p}$}{Fixed cost of technology $p$ \nomunit{\unit{\Currency\per\MW\per\year}}}
\nomenclature[C15]{$\invCost{p}$}{Investment cost of technology $p$ \nomunit{\unit{\Currency\per\MWh}}}
\nomenclature[C16]{$\transmissionCost$}{Investment cost for constructing transmission lines between regions \nomunit{\unit{\Currency\per\MWh\per\km}}}
\nomenclature[C17]{$\sheddingCost$}{Load shedding cost \nomunit{\unit{\Currency\per\MWh}}}
\nomenclature[C18]{$\crf{p}$}{Capital recovery factor for technology $p$}
\nomenclature[C19]{$\crfTransmission$}{Capital recovery factor of transmission lines $p$}
\nomenclature[C20]{$\dischargeTime$}{Battery full discharge time \nomunit{\unit{\hour}}}
\nomenclature[C21]{$\carbonTax$}{Carbon tax \nomunit{\unit{\Currency\per\TonCOtwo}}}
\nomenclature[C22]{$\prob{s}$}{Probability of scenario $s$}
\nomenclature[C23]{$\OperationalCost{s}$}{Operational cost under scenario $s$ \nomunit{\unit{\Currency\per\year}}}
\nomenclature[C24]{$\TotalEmission{s}$}{Annual \ch{CO2} emissions under scenario $s$ \nomunit{\unit{\TonCOtwo\per\year}}}
\nomenclature[C25]{$\TotalLoadShedding{s}$}{Annual shed load under scenario $s$  \nomunit{\unit{\MWh}}}
\nomenclature[C26]{$\FixedCost$}{Annual fixed cost   \nomunit{\unit{\Currency\per\year}}}
\nomenclature[C27]{$\InvestmentCost$}{Levelized investment on technologies  \nomunit{\unit{\Currency\per\year}}}
\nomenclature[C28]{$\SC$}{Levelized total costs  \nomunit{\unit{\Currency\per\year}}}

\nomenclature[D01]{$\Electricity{s}{r}{p}{t}$}{Electricity generated by technology $p$ in region $r$ under scenario $s$ at time step $t$  \nomunit{\unit{\MWh}}}
\nomenclature[D02]{$\Capacity{r}{p}$}{Capacity of technology $p$ in region $r$  \nomunit{\unit{\MW}}}
\nomenclature[D03]{$\TransmissionCapacity{r}{r'}$}{Transmission capacity between regions $r$ and $r'$  \nomunit{\unit{\MW}}}
\nomenclature[D04]{$\PlannedOutage{s}{r}{t}$}{Planned outage for nuclear plants in region $r$ under scenario $s$ at time step $t$  \nomunit{\unit{\MWh}}}
\nomenclature[D05]{$\ReservoirLevel{s}{r}{t}$}{Hydro reservoir level in region $r$ for scenario $s$ and time step $t$  \nomunit{\unit{\MWh}}} 
\nomenclature[D06]{$\Transmission{s}{r}{r'}{t}$}{Electricity transmitted from region $r$ to region $r'$ under scenario $s$ at time step $t$ \nomunit{\unit{\MWh}}}
\nomenclature[D07]{$\ElectricityShedding{s}{r}{t}$}{Shed load in region $r$ under scenario $s$ at time step $t$  \nomunit{\unit{\MWh}}}
\nomenclature[D08]{$\BatteryLevel{s}{r}{t}$}{Charge stored in the battery in region $r$ under scenario $s$ at time step $t$  \nomunit{\unit{\MWh}}}

\vspace{-4em}
\renewcommand{\nomname}{}
\printnomenclature

\section{\WU mathematical model} \label{sec:stochastic_model}

The mathematical model of the proposed capacity expansion problem with uncertainty in wind and solar power outputs (\WU model) is as follows.

\begin{align} 
    & \text{\Stochastic model:} \nonumber \\
    & \min \SC = \InvestmentCost + \FixedCost + \sum_{s \in \Scenarios} \prob{s} \left(\OperationalCost{s} + \carbonTax \cdot \TotalEmission{s} + \sheddingCost \cdot \TotalLoadShedding{s}\right) \label{eq:objective_function} \\
    &\text{Subject to:} \nonumber \\
    &\Capacity{r}{p} \le \MaxCapacity{r}{p}, \quad \forall r \in \Regions, \; p \in \Technologies \label{eq:max_capacity}\\
    & \sum_{p \in \Plants \cup \{\Inverter\}} \Electricity{s}{r}{p}{t} - \Electricity{s}{r}{\Battery}{t} + \sum_{\substack{r' \in \Regions \\ (r, r') \in \ValidLinks}} \left(\transmissionEfficiency{r'}{r} \cdot \Transmission{s}{r'}{r}{t} -  \Transmission{s}{r}{r'}{t}\right) \nonumber \\
     & \qquad \geq \timeStep \cdot \load{r}{t} - \ElectricityShedding{s}{r}{t}, \quad \forall s \in \Scenarios, \; r \in \Regions, \; t \in \TimeSteps \label{eq:electricity_demand} \\
    & \ElectricityShedding{s}{r}{t} \leq \sheddingCapRatio \cdot \timeStep \cdot \load{r}{t}, \quad \forall s \in \Scenarios, \; r \in \Regions, \; t \in \TimeSteps \label{eq:shed_load_capacity} \\
    & \sum_{t \in \TimeSteps} \PlannedOutage{s}{r}{t} = \plannedOutageRate \cdot \Card{\TimeSteps} \cdot \timeStep \cdot \Capacity{r}{\Nuclear}, \quad \forall s \in \Scenarios, \; r \in \Regions \label{eq:planned_outage_nuclear} \\
    & \ReservoirLevel{s}{r}{t} \leq \left\{\ReservoirLevel{s}{r}{\timeStep \cdot \Card{\TimeSteps}} \text{, if } h = 1;  \ReservoirLevel{s}{r}{h-\timeStep} \text{, otherwise} \right\} \nonumber \\
    & \qquad + \timeStep \cdot \resInflow{r}{t} - \Electricity{s}{r}{\Hydro}{t}, \quad \forall s \in \Scenarios, \; r \in \Regions, \; t \in \TimeSteps \label{eq:reservoir_balance} \\
    & \BatteryLevel{s}{r}{t} \leq \left\{\BatteryLevel{s}{r}{\timeStep \cdot \Card{\TimeSteps}} \text{, if } h = 1;  \BatteryLevel{s}{r}{h-\timeStep} \text{, otherwise} \right\} \nonumber \\
    & \qquad + \efficiency{\Inverter} \cdot \Electricity{s}{r}{\Battery}{t} - \frac{\Electricity{s}{r}{\Inverter}{t}}{\efficiency{\Inverter}}, \quad \forall s \in \Scenarios, \; r \in \Regions, \; t \in \TimeSteps \label{eq:electricity_storage_balance} \\
    & \Electricity{s}{r}{\Nuclear}{t} \leq \timeStep \cdot \Capacity{r}{\Nuclear} - \PlannedOutage{s}{r}{t}, \quad \forall s \in \Scenarios, \; r \in \Regions, \; t \in \TimeSteps \label{eq:gen_capacity_nuclear} \\
    & \Electricity{s}{r}{p}{t} \leq \timeStep \cdot \cf{s}{r}{p}{t} \cdot \Capacity{r}{p}, \quad \forall s \in \Scenarios, \; r \in \Regions, \nonumber \\
    & \qquad p \in \{\Wind, \PV\}, \; t \in \TimeSteps \label{eq:gen_capacity_wind_PV} \\
    & \Electricity{s}{r}{p}{t} \leq \timeStep \cdot \Capacity{r}{p}, \quad \forall s \in \Scenarios, \; r \in \Regions, \; p \in \{\Gas, \Hydro\}, \; t \in \TimeSteps \label{eq:gen_capacity_others} \\
    & \Electricity{s}{r}{\Inverter}{t} +  \Electricity{s}{r}{\Battery}{t} \leq \timeStep \cdot \Capacity{r}{\Inverter}, \quad \forall s \in \Scenarios, \; r \in \Regions, \; t \in \TimeSteps \label{eq:gen_capacity_inverter} \\
    & \dischargeTime \cdot \Capacity{r}{\Inverter} \leq \Capacity{r}{\Battery}, \quad \forall r \in \Regions, \; t \in \TimeSteps \label{eq:capacity_battery} \\
    & \Transmission{s}{r}{r'}{t} + \Transmission{s}{r'}{r}{t} \leq \timeStep \cdot \TransmissionCapacity{r}{r'}, \quad \forall s \in \Scenarios, \; r, r' \in \Regions, \nonumber \\ 
    & \qquad r' > r, \; (r, r') \in \ValidLinks, \; t \in \TimeSteps \label{eq:electricity_transmission_capacity} \\
    & \ReservoirLevel{s}{r}{t} \leq \hydroReservoir{r}, \quad \forall s \in \Scenarios, \; r \in \Regions, \; t \in \TimeSteps \label{eq:storage_reservoir} \\
    & \BatteryLevel{s}{r}{t} \leq \Capacity{r}{\Battery}, \quad \forall s \in \Scenarios, \; r \in \Regions, \; t \in \TimeSteps \label{eq:storage_battery} \\
    & \TotalLoadShedding{s} = \sum_{r \in \Regions}\sum_{t \in \TimeSteps} \ElectricityShedding{s}{r}{t}, \quad \forall s \in \Scenarios \label{eq:total_load_shedding} \\
    & \TotalEmission{s} = \sum_{r \in \Regions}\sum_{p \in \Technologies \setminus \left\{ \Battery \right\}}\sum_{t \in \TimeSteps}\emissionFactor{p} \cdot \frac{\Electricity{s}{r}{p}{t}}{\efficiency{p}}, \quad \forall s \in \Scenarios \label{eq:total_CO2_emissions} \\
    & \OperationalCost{s} = \sum_{r \in \Regions}\sum_{p \in \Technologies \setminus \left\{ \Battery \right\}}\sum_{t \in \TimeSteps} \left(\frac{\fuelCost{p}}{\efficiency{p}} + \varCost{p}\right)\Electricity{s}{r}{p}{t}, \quad \forall s \in \Scenarios \label{eq:operational_cost} \\
    & \FixedCost = \sum_{r \in \Regions}\sum_{p \in \Technologies} \fixedCost{p} \cdot  \Capacity{r}{p} \label{eq:fixed_cost} \\
    & \InvestmentCost = \sum_{r \in \Regions}\sum_{p \in \Technologies} \invCost{p} \cdot \crf{p} \cdot  \Capacity{r}{p} + \sum_{\substack{r, r' \in \Regions \\ r' > r}} \transmissionCost \cdot \crfTransmission \cdot \distance{r}{r'} \cdot \TransmissionCapacity{r}{r'} \label{eq:investment_cost} \\
    & \PlannedOutage{s}{r}{t}, \ReservoirLevel{s}{r}{t}, \Transmission{s}{r}{r'}{t}, \BatteryLevel{s}{r}{t} \geq 0, \quad \forall s \in \Scenarios, \; r, r' \in \Regions, \nonumber \\ 
    & \qquad (r, r') \in \ValidLinks, \; t \in \TimeSteps \label{eq:non_negativity_1} \\
    & \TransmissionCapacity{r}{r'} \geq 0, \quad \forall r, r' \in \Regions, \; r' > r, \; (r, r') \in \ValidLinks \label{eq:non_negativity_2} \\
    & \Capacity{r}{p} \geq 0, \quad \forall r \in \Regions, \; p \in \Technologies \label{eq:non_negativity_3} \\
    & \Electricity{s}{r}{p}{t} \geq 0, \quad \forall s \in \Scenarios, \; r \in \Regions, \; p \in \Technologies \setminus \left\{ \Battery \right\}, \; t \in \TimeSteps \label{eq:non_negativity_4}
\end{align}

The objective function~(\ref{eq:objective_function}) minimizes the levelized total cost, which includes the levelized investment cost along with the average operational, fixed, \ch{CO2} emission, and load shedding costs across different scenarios. Constraint~(\ref{eq:max_capacity}) limits the capacity of each technology in a region so that it does not exceed the maximum potential capacity of that region. Constraint~(\ref{eq:electricity_demand}) ensures that the total electricity supply from generation, storage, and transmission meets the demand after accounting for any shed load. Constraint~(\ref{eq:shed_load_capacity}) prevents the shed load from exceeding a certain fraction of the electricity demand. Constraint~(\ref{eq:planned_outage_nuclear}) enforces that the planned outage duration equals a certain fraction of the time in one year. Constraint~(\ref{eq:reservoir_balance}) states that the electrical energy stored in hydro reservoirs at the end of each time step is equal to that of the previous time step plus the water inflow (in energy terms) minus the electrical energy withdrawn from the reservoir. Similarly, constraint~(\ref{eq:electricity_storage_balance}) ensures that the amount of electrical energy stored in batteries at the end of each time step equals that of the previous time step plus the added charge minus the discharged amount. Note that both constraints~(\ref{eq:reservoir_balance}) and~(\ref{eq:electricity_storage_balance}) form a closed loop in the storage balance constraints by linking the last time step to the first time step. Constraint~(\ref{eq:gen_capacity_nuclear}) restricts the electricity produced by nuclear plants to be less than or equal to their capacity minus planned outages. Constraint~(\ref{eq:gen_capacity_wind_PV}) ensures that the wind and solar generation does not exceed their capacity, adjusted by the corresponding capacity factor for each time step. Constraint~(\ref{eq:gen_capacity_others})--(\ref{eq:capacity_battery}) ensures that the electricity generated by gas turbines, hydropower plants, the electricity used for battery charging, and the electricity discharged from batteries do not exceed their respective capacities. Constraint~(\ref{eq:electricity_transmission_capacity}) states that the electricity transmitted between regions does not exceed the transmission line capacity. Constraint~(\ref{eq:storage_reservoir}) restricts the stored electricity in the hydro reservoir to its maximum storage capacity in each region and constraint~(\ref{eq:storage_battery}) ensures that the electricity stored in the batteries does not exceed their storage capacity. Constraints~(\ref{eq:total_load_shedding})--(\ref{eq:investment_cost}) respectively compute the annual load shedding, \ch{CO2} emissions, operating and fixed costs, and levelized investment costs. Finally, the non-negativity requirement of the decision variables is stated in constraints~(\ref{eq:non_negativity_1})--(\ref{eq:non_negativity_4}).

\section{Mathematical model for simulating the solutions} \label{sec:simulation_model}

The following mathematical model is applied to evaluate a given capacity layout $\CapacityBarBold$ using simulation. 

\begin{align}
    &\hypertarget{SimulationProblem}{\text{Simulation problem}}(\noNuclearPlants, \CapacityBarBold, s)\text{:} \nonumber \\
    & \SC_s  = \InvestmentCost + \FixedCost + \min \left\{\OperationalCost{s} + \carbonTax \cdot \TotalEmission{s} + \sheddingCost \cdot \TotalLoadShedding{s} + \loadLossCost \cdot \LoL_{s} \right\} \label{eq:simulation_objective}\\ 
    & \text{subject to: (\ref{eq:shed_load_capacity})--(\ref{eq:electricity_storage_balance}), (\ref{eq:gen_capacity_wind_PV})--(\ref{eq:non_negativity_4}), (\ref{eq:gen_capacity_nuclear_robust}), (\ref{eq:non_negativity_5})}  \nonumber \\
    & \sum_{p \in \Plants \cup \{\Inverter\}} \Electricity{s}{r}{p}{t} - \Electricity{s}{r}{\Battery}{t} + \sum_{\substack{r' \in \Regions \\ (r, r') \in \ValidLinks}} \left(\transmissionEfficiency{r'}{r} \cdot \Transmission{s}{r'}{r}{t} -  \Transmission{s}{r}{r'}{t}\right) \nonumber \\ 
    & \qquad \geq \timeStep \cdot \load{r}{t} - \ElectricityShedding{s}{r}{t} - \ElectricityLoadLoss{s}{r}{t}, \quad \forall r \in \Regions, \; t \in \TimeSteps \label{eq:electricity_demand_simulation} \\
    & \sum_{r \in \Regions}{\UnPlannedOutage{s}{r}{t}} = \oneNuclearPlantSize \cdot \timeStep \cdot \sum_{i=1}^{\noNuclearPlants}\simulatedOutageSamples{s}{i}{t}, \quad \forall  t \in \TimeSteps \label{eq:unplanned_load_loss_simulation} \\
    & \Capacity{r}{p} = \CapacityBar{r}{p}, \quad \forall r \in \Regions, \; p \in \Technologies \label{eq:capacity_equality}\\
    & \LoL_{s} = \sum_{r \in \Regions}\sum_{t \in \TimeSteps}\ElectricityLoadLoss{s}{r}{t} \label{eq:loss_of_load}\\
    & \ElectricityLoadLoss{s}{r}{t} \ge 0, \quad  \forall r \in \Regions, \; t \in \TimeSteps \label{eq:non_negativity_6}
\end{align}

In the \hyperlink{SimulationProblem}{Simulation problem}, Objective function~\eqref{eq:simulation_objective} minimizes the total system cost, plus the penalty for LoL. Note that the investment and fixed costs ($\InvestmentCost$ and $\FixedCost$) are dependent only on the given capacity mix and can be pre-computed; this allows them to be excluded from the optimization. Constraint~\eqref{eq:electricity_demand_simulation} is analogous to Constraint~\eqref{eq:electricity_demand} but incorporates the loss of load. Constraint~\eqref{eq:unplanned_load_loss_simulation} calculates unplanned nuclear power outages based on the generated hourly power outages and the unit capacity of nuclear plants. Note that Algorithm~\ref{func:GenerateOutageSamples}($\noNuclearPlants$) is applied to generate (simulate) the time series for nuclear power outages (denoted by $\simulatedOutageSamples{s}{i}{t}$). Constraint~\eqref{eq:capacity_equality} ensures that the capacity mix remains fixed. Constraint~\eqref{eq:loss_of_load} computes the annual loss of load for a given simulation year, and Constraint~\eqref{eq:non_negativity_6} enforces non-negativity for the loss of load variable.

\section{Description of additional parameters in the case study} \label{sec:additional_parameters}

The investment cost, fixed cost, variable cost, fuel cost, efficiency, lifetime, and emission factor of each technology are given in Table~\ref{tab:cost_parameters}, which are mainly derived from the Danish Energy Agency \citep{DEA2024}. For the transmission grids, the investment cost, $\transmissionCost$, is \SI{0.4}{\Currency \per \KW \per \km} with a lifetime equal to 40 years \citep{hagspiel2014cost}. The transmission efficiency between the regions, $\transmissionEfficiency{r}{r'}$, is calculated by $(1- \lossPerKm)^{\distance{r}{r'} / 1000}$ where $\lossPerKm$ is the electricity loss per kilometer and is assumed to be 0.016. A discharge time of four hours is assumed for the batteries, i.e., $\dischargeTime = 4 \unit{\hour}$ \citep{cole2019cost}. 

The capital recovery factor of each technology is calculated by $\frac{\discountRate}{1 - (1 + \discountRate)^{-\lifeTime}}$, where $\discountRate$ is the discount rate and is assumed to be 0.05, and $\lifeTime$ is the life of the technology. We assume that the total shed load at each time step should not exceed 5\% of the electricity demand, i.e., $\sheddingCapRatio = 0.05$ \citep{van2019cost}. Additionally, the carbon tax $\carbonTax$, load shedding cost $\sheddingCost$, and penalty cost for loss of load (in the simulation model) $\loadLossCost$ are set to \SI{150}{\Currency \per \TonCOtwo}, \SI{1000}{\Currency \per \MWh}, and \SI{10000}{\Currency \per \MWh}, respectively. The planned outage rate, $\plannedOutageRate$, is also assumed to be 15\%. Each region has some requirements to ensure that investment in solar and wind technologies does not exceed certain limits, due to the limited land availability. Additionally, the hydropower and reservoir capacity are assumed to be fixed in each region due to environmental regulations. The limitations on solar and wind capacities, along with the existing capacity of the hydropower technology and reservoir in the regions, are given in Table~\ref{tab:max_capacity}; it should be mentioned that for the rest of the technologies, no capacity limitation is applied.

\section{Sensitivity analysis based on \Deterministic and \Stochastic models} \label{sec:sensitivity analysis}

\begin{figure}[htbp]
    \centering
        \includegraphics[width=0.5\columnwidth]{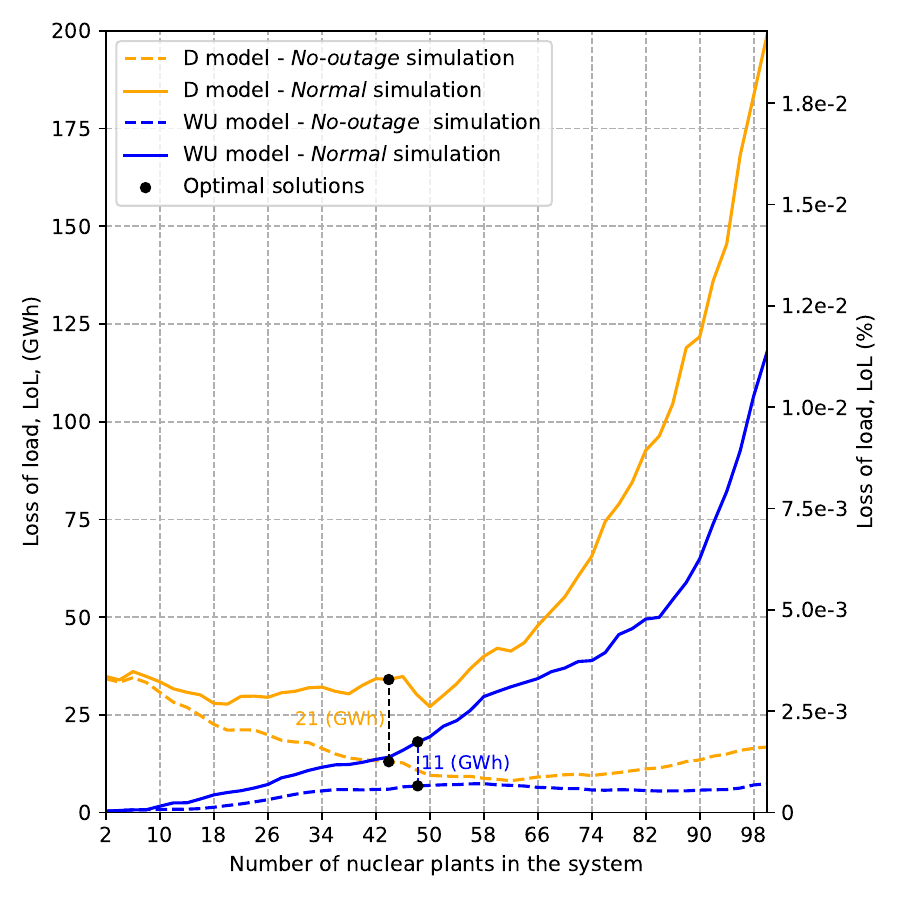}
    \caption{Loss of load (LoL) obtained from the simulation for the \D (\Deterministic) and \WU (\Stochastic) models plotted against different numbers of nuclear plants in the system.}
    \label{fig:LoL_plot_different_n}
\end{figure}

\begin{figure}[htbp]
    \centering
    \includegraphics[width=0.5\columnwidth]{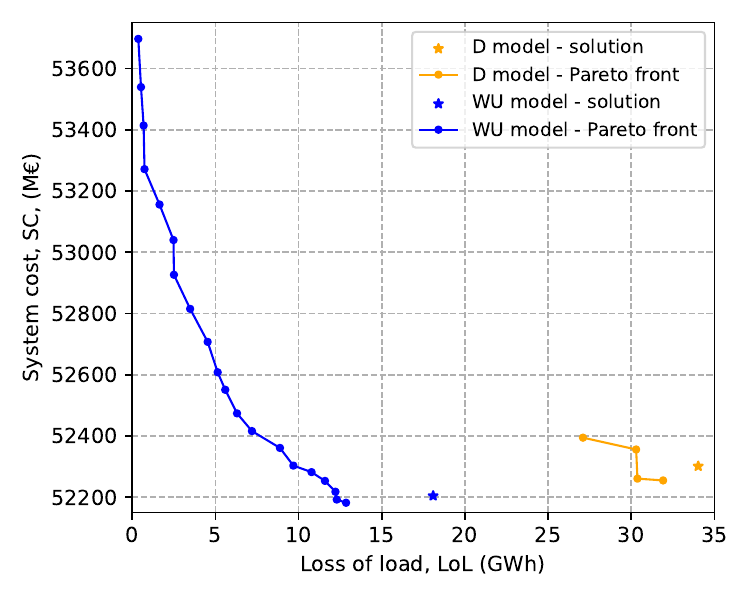}
    \caption{Pareto front of the \D (\Deterministic) and \WU (\Stochastic) models' solutions, illustrating the trade-off between the $\SC$ and LoL over a 41-year simulation.}
    \label{fig:pareto_front_stochastich_and_deterministic}
\end{figure}

For the sensitivity analysis, an additional constraint was added to the model presented in Section~\ref{sec:StochasticModel} to fix the total number of nuclear plants in the system. In Figure~\ref{fig:LoL_plot_different_n}, the dashed lines show the simulation results with only uncertainty in weather conditions, and the solid lines correspond to the simulation results with uncertainty in both weather conditions and nuclear power availability. As seen in this figure, increasing the penetration of nuclear power plants makes the system more vulnerable. Moreover, referring to Figure~\ref{fig:LoL_plot_different_n}, under the \NormalSimulation simulation, the LoL of the \Deterministic model's solution is always higher than that of the \Stochastic model's solution. Thus, the \Stochastic model provides some hedging also for nuclear power uncertainty (red line vs. blue line in Figure~\ref{fig:LoL_plot_different_n}), but the LoL still increases proportionally with increasing nuclear penetration. Based on the sensitivity analysis results, a Pareto front was derived for both models, as shown in Figure~\ref{fig:pareto_front_stochastich_and_deterministic}. This figure demonstrates that the \Stochastic model outperforms the \Deterministic model in terms of cost and LoL. 
\SetCommentSty{mycommfont}
\section{Algorithms}\label{sec:appendix_algorithms}

\begin{function}[h]
    \caption{SolveRobustProblem($\mathcal{A}, \noSamplesforPercentiles, \noSamplesinRobustModel$, $s$)} \label{func:SolveRobustProblem}
    \DontPrintSemicolon
    \SetArgSty{textrm}
    $\outageSamplesBold \gets \textbf{GenerateOutageSamples}(\noSamplesinRobustModel)$ \tcp{These outage samples are generated and used in the Master problem and Sub-problem.}
    $\text{pareto} \gets \emptyset$ \tcp{Pareto front is initialized to empty.}
    \tcc{For every confidence level ($\confidenceLevel$) and nuclear plant size ($\noNuclearPlants$), the robust model is solved.}
    \For{$(\noNuclearPlants \in \Sizes) \textbf{ and } (\confidenceLevel \in \mathcal{A})$}{
        $\SelectedSampleBar_{i, s} \gets 0, \quad \forall i \in \left\{1, 2, \ldots, \noSamplesinRobustModel \right\}, \; s \in \Scenarios$ \;
        Solve $\hyperlink{MasterProblem}{\text{Master problem}}(\noNuclearPlants, \SelectedSampleBarBold)$ \;
        $\CapacityBar{r}{p}  \gets \Capacity{r}{p}, \quad \forall r \in \Regions, \; p \in \Technologies$ \;
        $(\annualOutagePercentile, \maxNoOutagePercentile) \gets \protect\textbf{\ref{func:Percentiles}}(\noNuclearPlants, \confidenceLevel, \noSamplesforPercentiles)$ \tcp{Annual and maximum unplanned outages $(\annualOutagePercentile, \maxNoOutagePercentile)$ are calculated.}
        \tcc{For each scenario, the sub-problem is solved to identify the start time and duration of nuclear power outages that may result in a loss of load risk.}
        \For{$s \in \Scenarios$}{
            Solve $\hyperlink{SubProblem}{\text{Sub-problem}}(\noNuclearPlants, \annualOutagePercentile, \maxNoOutagePercentile, \noSamplesinRobustModel, \CapacityBarBold, s)$\;
            $\SelectedSampleBar_{i, s} \gets \SelectedSample_{i}, \quad \forall i \in \left\{1, 2, \ldots, \noSamplesinRobustModel \right\}$ \;
        }
        Solve $\hyperlink{MasterProblem}{\text{Master problem}}(\noNuclearPlants, \SelectedSampleBarBold)$ \;
        $\bar{c} \gets 0$ and $\bar{l} \gets 0$\;
        $\CapacityBar{r}{p}  \gets \Capacity{r}{p}, \quad \forall r \in \Regions, \; p \in \Technologies$ \;
        \tcc{Solution $\CapacityBarBold$ is simulated to accurately assess the current cost ($\bar{c}$) and current loss of load ($\bar{l}$).}
        \For{$s \in \ScenariosSimulation$}{
            $\simulatedOutageSamplesBold_{s} \gets \textbf{GenerateOutageSamples}(\noNuclearPlants)$ \;
            Solve $\hyperlink{SimulationProblem}{\text{Simulation problem}}(\noNuclearPlants, \CapacityBarBold)$ and obtain $\SC$\;
            $\bar{l} \gets \sum_{t \in \TimeSteps}\ElectricityLoadLoss{s}{r}{t}$ \;
            $\bar{c} \gets \bar{c} + \SC - \loadLossCost \cdot \sum_{t \in \TimeSteps}\ElectricityLoadLoss{s}{r}{t}$  \tcp{Loss of load is excluded from the cost.}
        }
        $\bar{c} \gets \frac{\bar{c}}{\Card{\ScenariosSimulation}}$ and  $\bar{l} \gets \frac{\bar{l}}{\Card{\ScenariosSimulation}}$\;
        
        \tcc{Update the Pareto front if solution $(\bar{c}, \bar{l}, \CapacityBarBold)$ dominates the current solutions in Pareto front}
        $\text{nonDominated} \gets \texttt{true}$\;
        \For{$\text{sol} \in \text{pareto}$}{
            \If{$(\bar{c} \geq \text{sol.cost} \; \textbf{and} \; \bar{l} \geq \text{sol.loss})$}{
                $\text{nonDominated} \gets \texttt{false}$\;
                \textbf{break}\;
            }
        }
        \If{$\text{nonDominated}$}{
            $\text{pareto} \gets \text{pareto} \setminus \left\{\text{sol} \in \text{pareto} : (\bar{c} \leq \text{sol.cost} \; \textbf{and} \; \bar{l} \leq \text{sol.loss})\right\}$\;
            $\text{pareto} \gets \text{pareto} \cup \left\{\text{solution} = (\bar{c}, \bar{l}, \CapacityBarBold)\right\}$\;    
        }
    }
    \Return $\text{pareto}$\;
\end{function}
\clearpage
\onecolumn
\section{Tables}\label{sec:appendix_tables}

\begin{table}[htbp]
    \centering
    \caption{Cost parameters, operational characteristics, and environmental impacts of the technologies considered in the case study. Cost assumptions are drawn mainly from the Danish Energy Agency \citep{DEA2024}.}
    \begin{threeparttable}
        \resizebox{\textwidth}{!}{
        \begin{tabularx}{\textwidth}{Xccccccc} \toprule
        & Investment Cost 
        & Fixed Cost   
        & Variable Cost 
        & Fuel Cost     
        & Efficiency 
        & Lifetime 
        & Emission Factor \\
        Plant      & \unit{\Currency \per \KW}      & \unit{\Currency \per \MW \per \year} & \unit{\Currency \per \MWh} & \unit{\Currency \per \MWh} & --         & \unit{\year} & \unit{\TonCOtwo \per \MWh}     \\ \midrule
        \Hydro     & 0           & 30000       & 0             & 0             & 0.9        & 80       & 0           \\
        \Gas       & 436         & 7893        & 4.79          & 32            & 0.43       & 25       & 0.202           \\
        \Wind      & 1090        & 15602       & 1.85          & 0             & 1          & 30       & 0               \\
        \PV        & 290         & 9900        & 0             & 0             & 1          & 40       & 0               \\
        \Battery   & 65\tnote{a} & 0           & 0             & 0             & 1          & 15       & 0               \\
        \Inverter  & 200         & 38000       & 0             & 0             & 0.92       & 15       & 0               \\
        \Nuclear   & 4000        & 126000      & 1.9           & 3             & 0.33       & 40       & 0               \\ \bottomrule 
        \end{tabularx}}
        \begin{tablenotes}
            \item[a] The unit for battery storage capacity is \unit{\Currency \per \KWh}.
            \item[b] To assess the impact of uncertainty in nuclear availability, we assume a relatively low cost for nuclear power to ensure it plays a significant role in the optimal energy mix. 
        \end{tablenotes}
    \end{threeparttable}
    \label{tab:cost_parameters}
\end{table}

\begin{table}[htbp]
    \centering
    \caption{Maximum wind and solar potential \citep{mattsson2021autopilot}, hydropower capacity, and reservoir storage capacity in each region \citep{entsoe2020}.} 
    \resizebox{\textwidth}{!}{
    \begin{tabularx}{\textwidth}{lXXXXXXX} \toprule
        & \multicolumn{7}{c}{Regions} \\ \cmidrule(lr){2-8}
        Technologies & \SEN & \SES & \DK & \DEN& \DES & \BNL& \PL \\ \midrule  
        Hydropower capacity, $\MaxCapacity{r}{\Hydro} (\unit{\GW})$ & 13.7 & 2.5 & 0 & 0.21 & 1.17 & 0 & 0.47 \\
        Wind capacity, $\MaxCapacity{r}{\Wind} (\unit{\GW})$ & 46 & 60 & 19 & 73 & 42 & 22 & 125 \\
        Solar capacity, $\MaxCapacity{r}{\PV}$ (\unit{\GW})& 11 & 39 & 59 & 233 & 142 & 59 & 432 \\
        Reservoir capacity, $\hydroReservoir{r} (\unit{\TWh})$ & 27.2 & 6.6 & 0 & 0.05 & 0.25 & 0 & 1.6 \\ \bottomrule 
    \end{tabularx}
    }
    \label{tab:max_capacity}
\end{table}

\begin{table}[htbp]
    \caption{Summary of the Pareto solutions obtained by solving the \WNU (\Robust) model.}
    \label{tab:pareto_solutions_summary}
    \resizebox{\textwidth}{!}{
        \begin{tabular}{lccccccccccccccc} \toprule
             & \multicolumn{8}{c}{Total capacity} & \multicolumn{5}{c}{System costs} & & \\ \cmidrule(lr){2-9} \cmidrule(lr){10-14} 
            \begin{tabular}[c]{@{}l@{}}Solution \\ \#\end{tabular} & 
            \begin{tabular}[c]{@{}c@{}}Hydropower\\ (\unit{\GW})\end{tabular} & 
            \begin{tabular}[c]{@{}c@{}}Gas\\ (\unit{\GW})\end{tabular} & 
            \begin{tabular}[c]{@{}c@{}}Wind\\ (\unit{\GW})\end{tabular} & 
            \begin{tabular}[c]{@{}c@{}}Solar\\ (\unit{\GW})\end{tabular} & 
            \begin{tabular}[c]{@{}c@{}}Nuclear\\ (\unit{\GW})\end{tabular} & 
            \begin{tabular}[c]{@{}c@{}}Transmission\\ (\unit{\GW})\end{tabular} & 
            \begin{tabular}[c]{@{}c@{}}Inverter\\ (\unit{\GW})\end{tabular} & 
            \begin{tabular}[c]{@{}c@{}}Battery\\ (\unit{\GWh})\end{tabular} & 
            \begin{tabular}[c]{@{}c@{}}Investment \\ (\unit{\MCurrency})\end{tabular} & 
            \begin{tabular}[c]{@{}c@{}}Fixed \\ (\unit{\MCurrency})\end{tabular} & 
            \begin{tabular}[c]{@{}c@{}}Operational \\ (\unit{\MCurrency})\end{tabular} & 
            \begin{tabular}[c]{@{}c@{}}Emissions \\ (\unit{\MCurrency})\end{tabular} & 
            \begin{tabular}[c]{@{}c@{}}Load shedding \\ (\unit{\MCurrency})\end{tabular} & 
            \begin{tabular}[c]{@{}c@{}}$\SC$ \\ (\unit{\MCurrency})\end{tabular} & 
            \begin{tabular}[c]{@{}c@{}}LoL \\ (\unit{\GWh})\end{tabular} \\ \midrule
                1 & 18 & 60 & 182 & 223 & 38 & 90 & 32 & 297 & 31028 & 12077 & 6065 & 2979 & 38 & 52186 & 6.32 \\ 
                2 & 18 & 59 & 177 & 219 & 40 & 89 & 31 & 279 & 30931 & 12165 & 6130 & 2925 & 39 & 52190 & 6.26 \\ 
                3 & 18 & 59 & 177 & 217 & 40 & 89 & 31 & 279 & 30867 & 12145 & 6180 & 2960 & 39 & 52191 & 6.25 \\ 
                4 & 18 & 60 & 182 & 224 & 38 & 92 & 32 & 294 & 31059 & 12080 & 6051 & 2970 & 36 & 52196 & 5.99 \\ 
                5 & 18 & 60 & 183 & 224 & 38 & 90 & 32 & 281 & 31045 & 12110 & 6045 & 2976 & 34 & 52209 & 5.59 \\ 
                6 & 18 & 60 & 178 & 218 & 40 & 88 & 32 & 274 & 30928 & 12183 & 6133 & 2936 & 33 & 52212 & 4.83 \\ 
                7 & 18 & 60 & 178 & 219 & 40 & 88 & 31 & 273 & 30976 & 12196 & 6097 & 2912 & 32 & 52213 & 4.76 \\ 
                8 & 18 & 60 & 178 & 218 & 40 & 88 & 31 & 270 & 30931 & 12177 & 6136 & 2941 & 32 & 52217 & 4.75 \\ 
                9 & 18 & 58 & 169 & 208 & 44 & 85 & 30 & 251 & 30799 & 12385 & 6223 & 2789 & 37 & 52232 & 4.52 \\ 
                10 & 18 & 60 & 178 & 218 & 40 & 88 & 32 & 268 & 30953 & 12203 & 6118 & 2930 & 29 & 52233 & 3.93 \\ 
                11 & 18 & 61 & 183 & 224 & 38 & 91 & 33 & 281 & 31130 & 12134 & 6006 & 2949 & 27 & 52245 & 3.73 \\ 
                12 & 18 & 61 & 183 & 224 & 38 & 91 & 33 & 290 & 31147 & 12134 & 5999 & 2938 & 26 & 52245 & 3.61 \\ 
                13 & 18 & 61 & 179 & 218 & 40 & 88 & 31 & 273 & 31064 & 12205 & 6062 & 2890 & 24 & 52246 & 2.82 \\ 
                14 & 18 & 60 & 171 & 206 & 44 & 85 & 29 & 240 & 30846 & 12374 & 6223 & 2799 & 27 & 52270 & 2.60 \\ 
                15 & 18 & 59 & 171 & 211 & 44 & 84 & 31 & 246 & 30972 & 12477 & 6090 & 2711 & 27 & 52276 & 2.35 \\ 
                16 & 18 & 62 & 180 & 219 & 40 & 90 & 31 & 269 & 31136 & 12227 & 6037 & 2867 & 21 & 52288 & 1.78 \\ 
                17 & 18 & 62 & 181 & 219 & 40 & 90 & 31 & 258 & 31174 & 12236 & 6007 & 2860 & 20 & 52298 & 1.45 \\ 
                18 & 18 & 61 & 172 & 207 & 44 & 84 & 29 & 240 & 30972 & 12390 & 6162 & 2765 & 19 & 52309 & 1.25 \\ 
                19 & 18 & 63 & 185 & 223 & 38 & 93 & 32 & 273 & 31257 & 12132 & 5974 & 2928 & 18 & 52310 & 1.23 \\ 
                20 & 18 & 61 & 176 & 215 & 42 & 88 & 31 & 253 & 31143 & 12383 & 6018 & 2764 & 21 & 52329 & 0.96 \\ 
                21 & 18 & 60 & 167 & 205 & 46 & 81 & 29 & 232 & 30935 & 12558 & 6169 & 2663 & 19 & 52344 & 0.75 \\ 
                22 & 18 & 61 & 169 & 203 & 46 & 82 & 28 & 222 & 30984 & 12531 & 6165 & 2670 & 16 & 52365 & 0.73 \\ 
                23 & 18 & 63 & 173 & 208 & 44 & 84 & 28 & 229 & 31054 & 12400 & 6143 & 2763 & 13 & 52373 & 0.69 \\ 
                24 & 18 & 65 & 184 & 217 & 40 & 91 & 29 & 248 & 31344 & 12228 & 5978 & 2843 & 11 & 52405 & 0.28 \\ 
                25 & 18 & 66 & 191 & 228 & 36 & 96 & 32 & 282 & 31523 & 12067 & 5865 & 2949 & 11 & 52415 & 0.27 \\ 
                26 & 18 & 66 & 188 & 222 & 38 & 94 & 30 & 257 & 31435 & 12130 & 5939 & 2913 & 10 & 52427 & 0.27 \\ 
                27 & 18 & 61 & 165 & 200 & 48 & 80 & 27 & 214 & 31049 & 12659 & 6149 & 2559 & 12 & 52427 & 0.24 \\ 
                28 & 18 & 65 & 180 & 212 & 42 & 87 & 29 & 237 & 31324 & 12352 & 5995 & 2764 & 9 & 52444 & 0.02 \\ 
                29 & 18 & 63 & 166 & 201 & 48 & 79 & 28 & 215 & 31210 & 12724 & 6053 & 2497 & 7 & 52490 & 0.00 \\ \bottomrule
        \end{tabular}
    }
\end{table}

\begin{table}[htbp]
    \caption{41 years simulation results considering unplanned nuclear power outages for the solutions depicted in Figures~\ref{fig:D_and_WU_map} and \ref{fig:robust_solutions}.} \label{tab:detailed_simulation_results}
    \resizebox{\textwidth}{!}{
        \begin{tabular}{lcccccccccccc} \toprule
            \multicolumn{1}{c}{} & \multicolumn{3}{l}{Deterministic solution}     & 
            \multicolumn{3}{c}{Stochastic solution}    & 
            \multicolumn{3}{c}{\begin{tabular}[c]{@{}c@{}} Robust solution\\ (Pareto solution \#1)\end{tabular}} & 
            \multicolumn{3}{c}{\begin{tabular}[c]{@{}c@{}} Robust solution\\ (Pareto solution \#29)\end{tabular}} \\ 
            \cmidrule(rl){2-4} \cmidrule(rl){5-7} \cmidrule(rl){8-10} \cmidrule(rl){11-13}
            Year                 
            & $\SC$ (\unit{\MCurrency}) 
            & LoL (\unit{\GWh})
            & LoL  (\%) 
            & $\SC$ (\unit{\MCurrency}) 
            & LoL (\unit{\GWh}) 
            & LoL  (\%) 
            & $\SC$ (\unit{\MCurrency})        
            & LoL (\unit{\GWh})      
            & LoL  (\%)       
            & $\SC$ (\unit{\MCurrency})        
            & LoL (\unit{\GWh})       
            & LoL  (\%)    \\ \midrule
                1979 & 51767 & 0 & 0 & 51639 & 0 & 0 & 51604 & 0 & 0 & 51962 & 0 & 0 \\
                1980 & 52154 & 0 & 0 & 52144 & 0 & 0 & 52301 & 0 & 0 & 52521 & 0 & 0 \\
                1981 & 52813 & 0 & 0 & 52711 & 1 &  $9.1\times 10^{-3}$  & 52644 & 0 & 0 & 52855 & 0 & 0 \\
                1982 & 53453 & 35 &  $3.3 \times 10^{-3}$  & 52731 & 39 &  $3.7\times 10^{-3}$  & 52940 & 0 & 0 & 52931 & 0 & 0 \\
                1983 & 49930 & 0 & 0 & 50081 & 0 & 0 & 49930 & 0 & 0 & 50534 & 0 & 0 \\
                1984 & 53102 & 0 & 0 & 52933 & 0 & 0 & 53099 & 0 & 0 & 53289 & 0 & 0 \\
                1985 & 53937 & 42 &  $4.0\times 10^{-3}$   & 53784 & 52 &  $5.0\times 10^{-3}$  & 53893 & 0 & 0 & 53904 & 0 & 0 \\
                1986 & 50395 & 0 & 0 & 50471 & 0 & 0 & 50358 & 0 & 0 & 50917 & 0 & 0 \\
                1987 & 54900 & 0 & 0 & 54393 & 0 & 0 & 54700 & 0 & 0 & 54678 & 0 & 0 \\
                1988 & 50113 & 0 & 0 & 50279 & 0 & 0 & 50005 & 0 & 0 & 50763 & 0 & 0 \\
                1989 & 50862 & 13 &  $1.3\times 10^{-3}$   & 50966 & 20 &  $1.9\times 10^{-3}$  & 50837 & 0 & 0 & 51292 & 0 & 0 \\
                1990 & 50527 & 0 & 0 & 50588 & 0 & 0 & 50597 & 0 & 0 & 51064 & 0 & 0 \\
                1991 & 52403 & 97 &  $9.3\times 10^{-3}$   & 52233 & 82 &  $7.9\times 10^{-3}$  & 52251 & 40 &  $3.8\times 10^{-3}$  & 52520 & 0 & 0 \\
                1992 & 51356 & 136 &  $1.3\times 10^{-3}$   & 51405 & 9 &  $8.2\times 10^{-3}$  & 50991 & 0 & 0 & 51600 & 0 & 0 \\
                1993 & 51409 & 20 &  $1.9\times 10^{-3}$   & 51499 & 0 & 0 & 51508 & 0 & 0 & 51795 & 0 & 0 \\
                1994 & 52129 & 0 & 0 & 51920 & 0 & 0 & 51850 & 0 & 0 & 52248 & 0 & 0 \\
                1995 & 53286 & 94 &  $9.0\times 10^{-3}$   & 52952 & 37 &  $3.5\times 10^{-3}$  & 52773 & 0 & 0 & 53132 & 0 & 0 \\
                1996 & 54535 & 273 &  $2.6\times 10^{-3}$   & 54130 & 156 &  $1.5\times 10^{-3}$  & 54328 & 122 &  $1.2\times 10^{-3}$  & 54248 & 0 & 0 \\
                1997 & 53429 & 0 & 0 & 53290 & 0 & 0 & 53435 & 0 & 0 & 53544 & 0 & 0 \\
                1998 & 51058 & 36 &  $3.5\times 10^{-3}$   & 51197 & 36 &  $3.5\times 10^{-3}$  & 51192 & 0 & 0 & 51630 & 0 & 0 \\
                1999 & 51367 & 0 & 0 & 51358 & 0 & 0 & 51208 & 0 & 0 & 51717 & 0 & 0 \\
                2000 & 51790 & 0 & 0 & 51702 & 0 & 0 & 51825 & 0 & 0 & 52055 & 0 & 0 \\
                2001 & 54360 & 152 &  $1.5\times 10^{-3}$   & 53701 & 100 &  $9.7\times 10^{-3}$  & 53758 & 40 &  $3.8\times 10^{-3}$  & 53826 & 0 & 0 \\
                2002 & 52735 & 0 & 0 & 52497 & 0 & 0 & 52453 & 0 & 0 & 52722 & 0 & 0 \\
                2003 & 52664 & 0 & 0 & 52559 & 0 & 0 & 52604 & 0 & 0 & 52794 & 0 & 0 \\
                2004 & 52207 & 0 & 0 & 52221 & 0 & 0 & 52062 & 0 & 0 & 52522 & 0 & 0 \\
                2005 & 51989 & 269 &  $2.6\times 10^{-3}$   & 51955 & 128 &  $1.2\times 10^{-3}$  & 52002 & 57 &  $5.5\times 10^{-3}$  & 52213 & 0 & 0 \\
                2006 & 51779 & 0 & 0 & 51784 & 0 & 0 & 51777 & 0 & 0 & 52050 & 0 & 0 \\
                2007 & 51330 & 0 & 0 & 51458 & 0 & 0 & 51342 & 0 & 0 & 51924 & 0 & 0 \\
                2008 & 51903 & 0 & 0 & 51830 & 0 & 0 & 51686 & 0 & 0 & 52182 & 0 & 0 \\
                2009 & 54971 & 0 & 0 & 54516 & 0 & 0 & 54714 & 0 & 0 & 54632 & 0 & 0 \\
                2010 & 53712 & 0 &  $1.3\times 10^{-3}$   & 53376 & 11 &  $1.0\times 10^{-3}$  & 53308 & 0 & 0 & 53504 & 0 & 0 \\
                2011 & 51450 & 0 & 0 & 51778 & 0 & 0 & 51720 & 0 & 0 & 52135 & 0 & 0 \\
                2012 & 52395 & 28 &  $2.7\times 10^{-3}$   & 52512 & 29 &  $2.7\times 10^{-3}$  & 52501 & 0 & 0 & 52776 & 0 & 0 \\
                2013 & 52976 & 0 & 0 & 52779 & 0 & 0 & 52875 & 0 & 0 & 53122 & 0 & 0 \\
                2014 & 53501 & 202 &  $1.9\times 10^{-3}$   & 53423 & 45 &  $4.3\times 10^{-3}$  & 53173 & 0 & 0 & 53339 & 0 & 0 \\
                2015 & 50633 & 0 & 0 & 50687 & 0 & 0 & 50756 & 0 & 0 & 51211 & 0 & 0 \\
                2016 & 53450 & 0 & 0 & 53244 & 0 & 0 & 53100 & 0 & 0 & 53333 & 0 & 0 \\
                2017 & 52576 & 0 & 0 & 52469 & 0 & 0 & 52204 & 0 & 0 & 52714 & 0 & 0 \\
                2018 & 52259 & 0 & 0 & 52240 & 0 & 0 & 52509 & 0 & 0 & 52563 & 0 & 0 \\
                2019 & 50754 & 0 & 0 & 50964 & 0 & 0 & 50822 & 0 & 0 & 51336 & 0 & 0 \\ \midrule
                Average  & 52301 & 34 &  $3.3\times 10^{-3}$   & 52205 & 18 &  $1.7\times 10^{-3}$  & 52186 & 6 &  $6.1\times 10^{-3}$  & 52490 & 0 & 0 \\
                Std.     & 1297 & 72 &  $6.9\times 10^{-3}$   & 1117 & 37 &  $3.5\times 10^{-3}$  & 1199 & 22 &  $2.1\times 10^{-3}$  & 1019 & 0 & 0 \\
                Max.     & 54971 & 273 &  $2.6\times 10^{-3}$   & 54516 & 156 &  $1.5\times 10^{-3}$  & 54714 & 122 &  $1.2\times 10^{-3}$  & 54678 & 0 & 0 \\ \bottomrule      
        \end{tabular}
    }
\end{table}
\clearpage
\onecolumn
\section{Figures}\label{sec:appendix_figures}

\begin{figure*}[htbp]
    \centering
    \begin{subfigure}[b]{0.45\linewidth}
        \includegraphics[width=\linewidth]{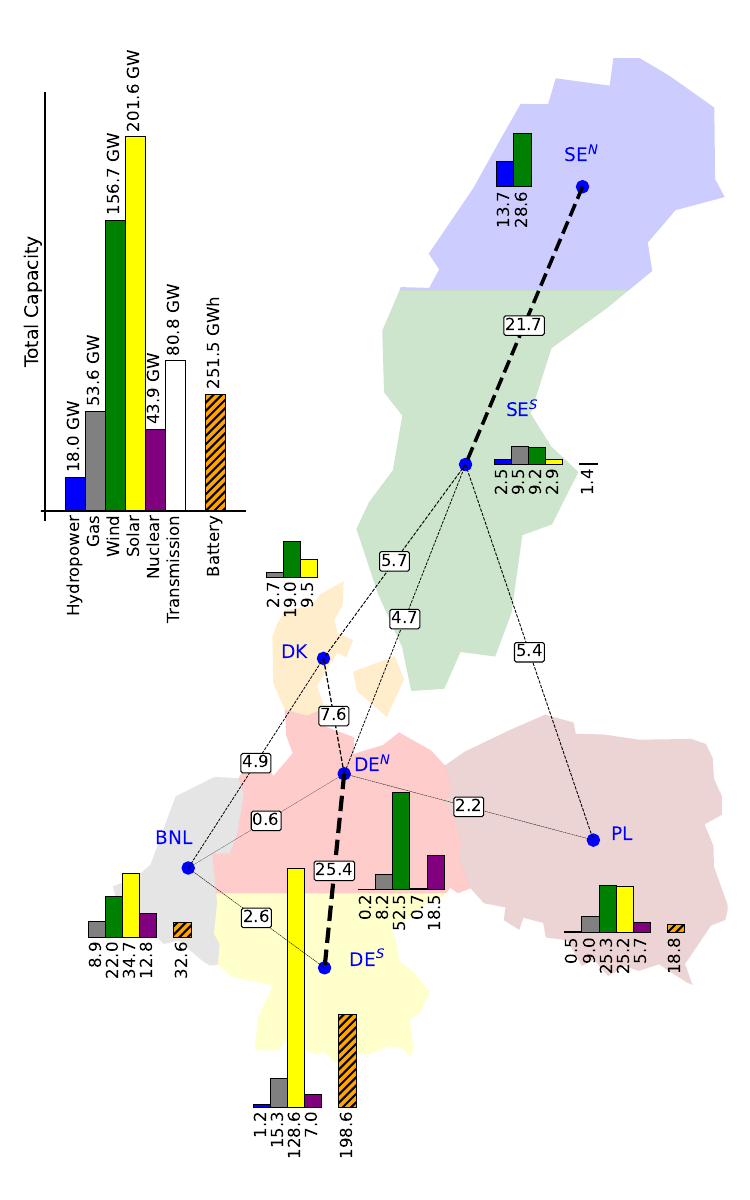}
        \caption{Deterministic solution, with $\SC$ = 52301~\unit{\MCurrency} and LoL = 34~\unit{\GWh}.} 
        \label{fig:D_map}
    \end{subfigure}
    \begin{subfigure}[b]{0.45\linewidth}
        \includegraphics[width=\linewidth]{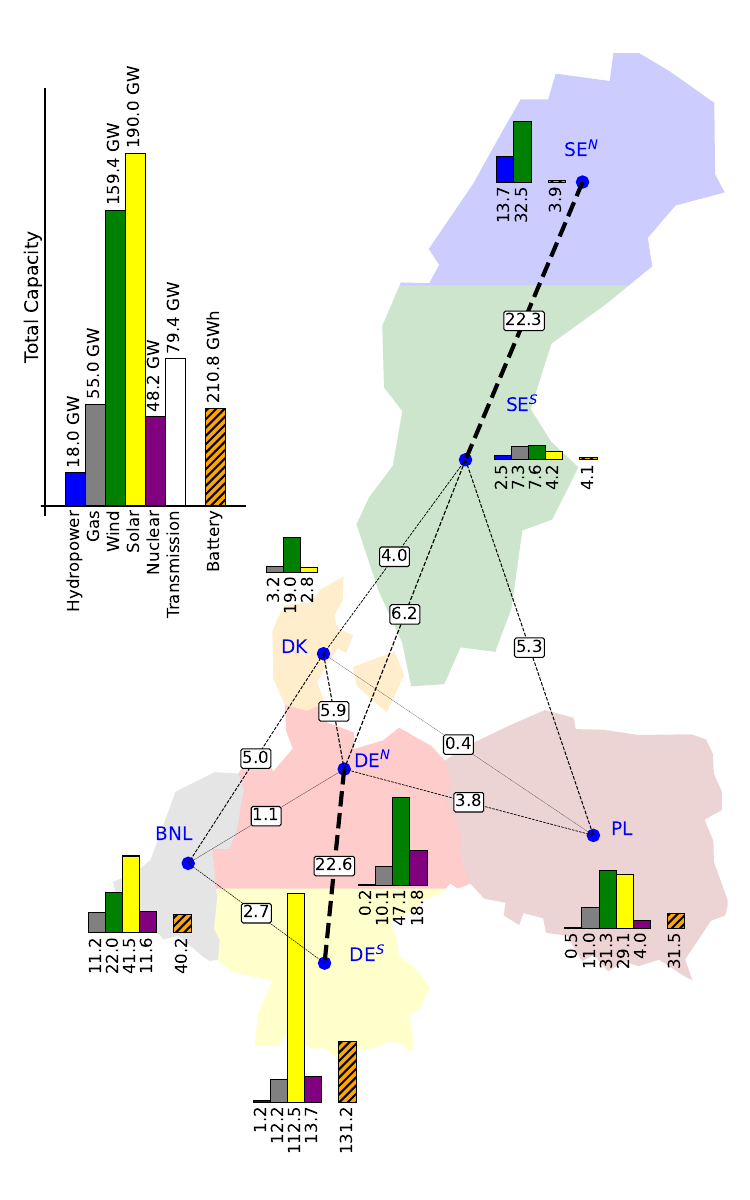}
        \caption{Stochastic solution, with $\SC$ = 52205~\unit{\MCurrency} and LoL = 18~\unit{\GWh}.} 
        \label{fig:WU_map}
    \end{subfigure}
    \caption{Optimal solution obtained by solving deterministic and stochastic models.}
    \label{fig:D_and_WU_map}
\end{figure*}

\begin{figure*}[htbp]
    \centering
    \begin{subfigure}[b]{0.45\linewidth}
        \includegraphics[width=\linewidth]{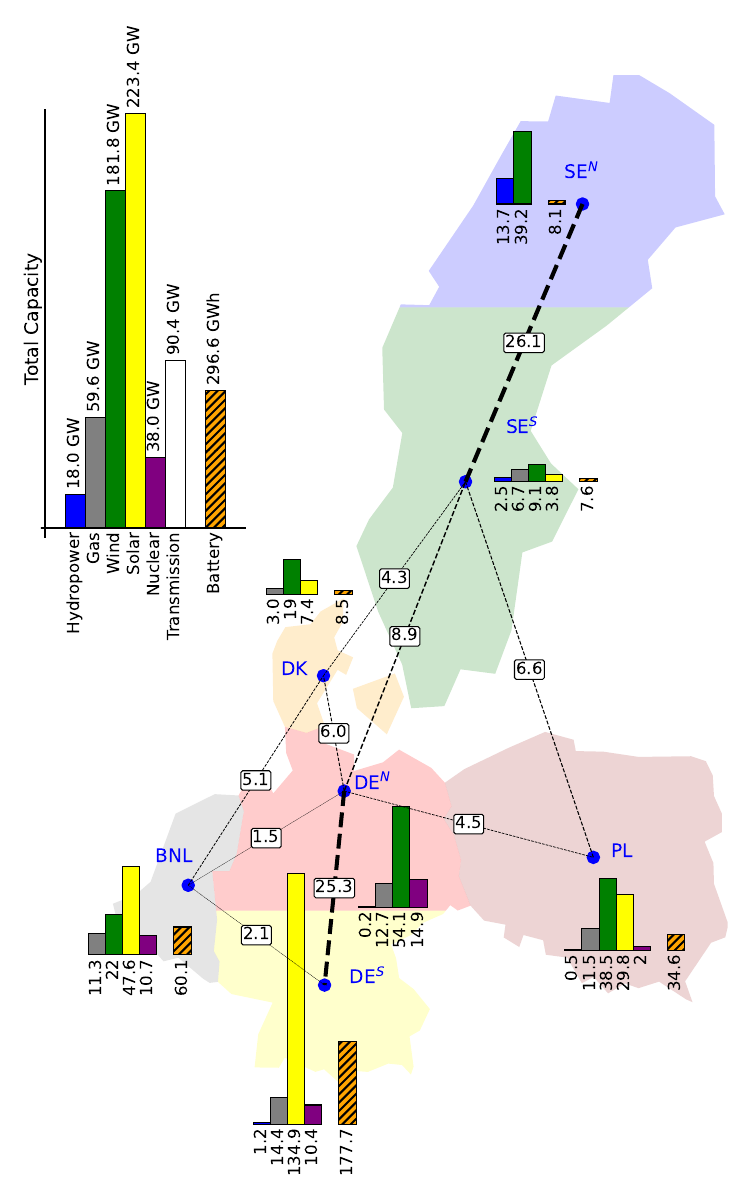}
        \caption{Pareto solution \#1, with $\SC$ = 52186~\unit{\MCurrency} and LoL = 6~\unit{\GWh}.} 
        \label{fig:robust_map_1}
    \end{subfigure}
    \begin{subfigure}[b]{0.45\linewidth}
        \includegraphics[width=\linewidth]{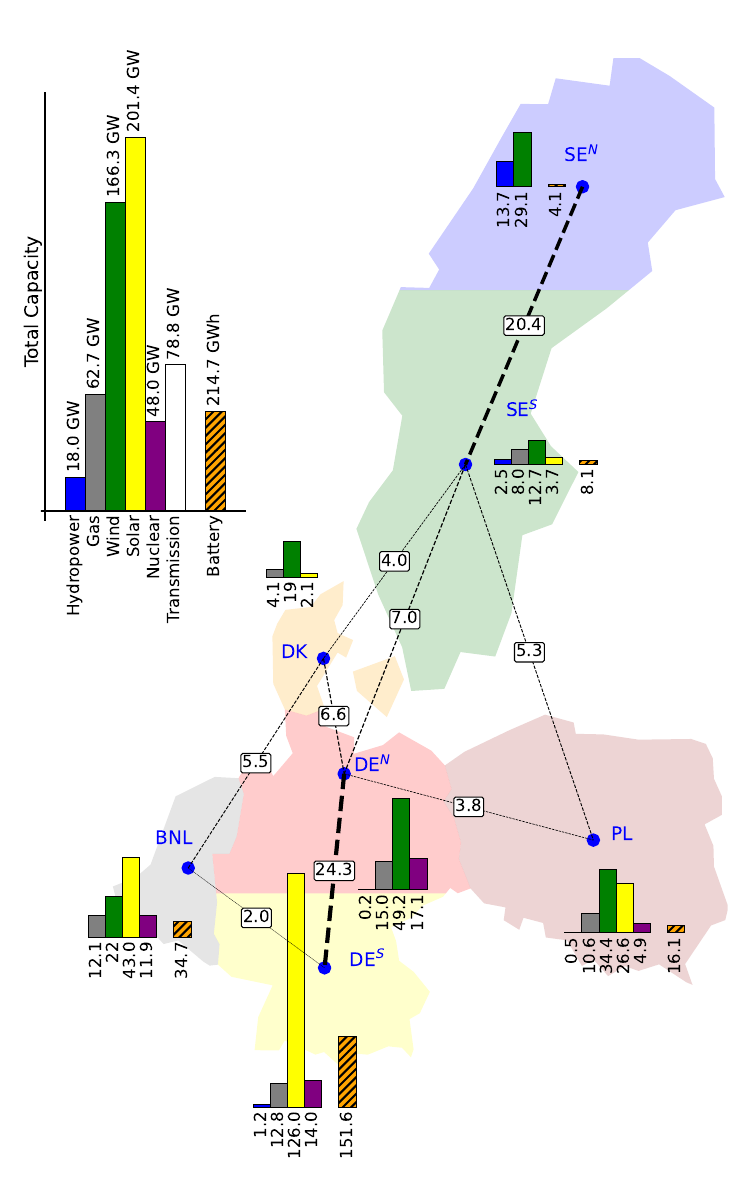}
        \caption{Pareto solution \#29, with $\SC$ = 52490~\unit{\MCurrency} and LoL = 0 ~\unit{\GWh}.} 
        \label{fig:robust_map_29}
    \end{subfigure}
    \caption{Two of the Pareto solutions obtained by the heuristic algorithm.}
    \label{fig:robust_solutions}
\end{figure*}

\end{document}